\documentclass[a4paper,12pt]{amsart}
\usepackage{amsmath,amssymb,mathrsfs,epsf}
\usepackage{hyperref}
\usepackage{fullpage}

\newtheorem{theo}{Theorem}[section]
\newtheorem{lemma}[theo]{Lemma}
\newtheorem{coro}[theo]{Corollary}
\newtheorem{prop}[theo]{Proposition}

\theoremstyle{definition}
\newtheorem{defi}[theo]{Definition}
\theoremstyle{remark}
\newtheorem{rem}[theo]{Remark}

\def\<#1,#2>{\langle #1,#2\rangle}
\def\braket#1{\mathinner{\langle{#1}\rangle}}

\renewcommand{\geq}{\geqslant}
\renewcommand{\leq}{\leqslant}

\newcommand{\R}{\mathbb{R}}
\newcommand{\E}{\mathbb{E}}
\newcommand{\N}{\mathcal{N}}
\newcommand{\pP}{\mathcal{P}}

\newcommand{\C}{\mathcal{C}}
\newcommand{\D}{\mathcal{D}}
\newcommand{\cX}{\mathcal{X}}
\newcommand{\sG}{\mathcal{G}}
\newcommand{\cU}{\mathcal{U}}

\newcommand{\Lip}{\operatorname{Lip}}
\newcommand{\mylip}[1]{\Lip (#1)}
\newcommand{\cS}{\mathcal{S}}
\newcommand{\sD}{\mathcal{D}}
\newcommand{\sym}{\operatorname{S}}
\newcommand{\symn}{\sym^n}
\newcommand{\symdim}[1]{\sym^{#1}}

\newcommand{\symnpos}{\symn_+}
\newcommand{\symndefpos}{\symn_{++}}

\newcommand{\Real}{\operatorname{Re}}

\newcommand{\cD}{\operatorname{D}}

\newcommand{\Tr}{\operatorname{Tr}}
\newcommand{\diag}{\operatorname{diag}}

\newcommand{\End}{\operatorname{End}}
\newcommand{\ssym}{\operatorname{Sym}}
\newcommand{\firstdef}[1]{{\em #1}}

\title[The contraction rate in Thompson metric of order-preserving flows]{The contraction rate in Thompson metric of order-preserving flows on a cone - application to generalized Riccati equations}

\author{St\'ephane Gaubert}
\address{INRIA and CMAP\\
         \'Ecole Polytechnique\\
          91128 Palaiseau C\'edex, France}
\email[]{Stephane.Gaubert@inria.fr}
\thanks{The authors were partially supported by the Arpege program of the French National Agency of Research (ANR), project ``ASOPT'', number ANR-08-SEGI-005 , by the Digiteo project DIM08 ``PASO'' number 3389.}
\author{Zheng QU}
\address{CMAP, INRIA and Fudan University\\
         \'Ecole Polytechnique\\
          91128 Palaiseau C\'edex, France}
\email[]{zheng.qu@polytechnique.edu}

\begin{document}
\begin{abstract}
We give a formula for the Lipschitz constant in Thompson's part metric of any order-preserving flow on the interior of a (possibly infinite dimensional) closed convex pointed cone. This provides an explicit form of a characterization of Nussbaum concerning non order-preserving flows.
As an application of this formula, we show that the flow of the generalized Riccati equation arising in stochastic linear quadratic control is a local contraction on the cone of positive definite matrices and characterize its Lipschitz constant by a matrix inequality. 
We also show that the same flow is no longer a contraction in other natural
Finsler metrics on this cone, including the standard invariant Riemannian metric.
This is motivated by a series of contraction properties
concerning the
standard Riccati equation, established by Bougerol,
Liverani, Wojtowski, Lawson, Lee and Lim:
we show that some of these properties do, and that some other do not, carry over
to the generalized Riccati equation.
\end{abstract}
\maketitle

\section{Introduction}
The standard discrete or differential Riccati equation arising in linear-quadratic control or optimal filtering problems has remarkable properties. 
In particular, Bougerol~\cite{MR1227540} proved that the 
standard discrete Riccati operator is non-expansive
in the invariant Riemannian metric on the set of positive definite matrices,
and that it is a strict contraction under controllability/observability
conditions.
%
%
Liverani and Wojtowski~\cite{liveraniw} proved that analogous contraction
properties hold with respect to Thompson's part metric.
These results, which were
obtained from algebraic properties of the linear symplectic semigroup associated
to a Riccati equation, are reminiscent
of Birkhoff's theorem in Perron-Frobenius theory (on the contraction
of positive linear operators sending a cone to its interior~\cite{birkhoff57}).
Lawson and Lim~\cite{MR2338433} generalized these results to the infinite
dimensional setting, and derived analogous contraction properties
for the flow of the differential Riccati equation
\begin{align}
\dot{P}= A'P+PA- P\Sigma P + Q, \qquad P(0)=G \enspace,
\label{e-standard}
\end{align}
where $A$ is a square matrix, $\Sigma,Q$ are positive
semi-definite matrices, and $G$ is a positive definite matrix.
Moreover, Lee and Lim~\cite{MR2399829} showed that the same contraction
properties hold more generally
for a family of Finsler metrics invariant under the action of the linear group (the latter metrics arise from
symmetric gauge functions). 

It is natural to ask whether contraction properties remain valid for more general
equations, like the following
constrained differential Riccati equation,
\begin{equation}\label{e-sriccnew}
\begin{array}{l}
\dot P=A'P+PA+C'PC+Q\\ \qquad \quad-(PB+C'PD+L')(R+D'PD)^{-1}(B'P+D'PC+L),\enspace 
\\
P(0)=G\\
R+D'PD \quad \text{positive definite,}
\end{array}
\end{equation}
which has received a considerable attention
in stochastic linear quadratic optimal control. The equation~\eqref{e-sriccnew} is known
as the {\em generalized Riccati differential equation} (GRDE)
or as the {\em stochastic Riccati differential equation}.
Up to a reversal of time,
it is a special case of the Backward stochastic Riccati
differential equation, which have extensively studied, see 
in particular~\cite{YongZhoubook99,ChengLiZhou98,aitrami01}.
The reader is referred specially to the monograph by Yong and Zhou~\cite{YongZhoubook99} for an introduction. 
Even for the simpler Riccati equation~\eqref{e-standard},
contractions
properties have not been established when the matrices $Q,\Sigma$ 
are not positive semi-definite, whereas this situation does occur in applications~\cite{curseofdim}.

In this paper, motivated by the analysis of the generalized Riccati equation,
we study the general question of computing the contraction
rate in Thompson's metric of an 
arbitrary order-preserving (time-dependent)
flow defined on a subset
of the interior of a closed convex and pointed cone in
a possibly infinite dimensional Banach space. 
Recall that the order associated with such a cone $\C$
is defined by $x\leq y\Leftrightarrow y-x\in \C$, and that the
Thompson metric can be defined on the interior of 
$\C$ by the formula
\[
d_T(x,y):=\log(\max \{M(x/y),M(y/x)\})
\]
where 
\[
M(x/y):=\inf\{t\in \R:ty\geq x\} =\sup_{\psi\in \C^*} \frac{\psi(x)}{\psi(y)}\enspace ,
\]
and $\C^*$ denotes the dual cone of $\C$.
More background can be found in~\S\ref{subsec-def-th}.

Our first main result 
can be stated as follows.
\begin{theo}\label{th-1}
Assume that the flow
of the differential equation $\dot{x}(t)=\phi(t,x(t))$ 
is order preserving with respect to the cone $\C$,
and let $\cU$ denote an open domain 
included in the interior of this cone such that $\lambda \cU \subset \cU$ holds for all $\lambda \in (0,1]$.
Then, the contraction rate of the flow over a time interval $J$, 
on the domain $\cU$, with respect to Thompson metric, is given by the formula
\begin{align}
\alpha:=-\sup_{s\in J,\; x\in \cU}  M\big((D\phi_s(x)x-\phi(s,x))/x\big) \enspace .
\label{e-explicit}
\end{align}
\end{theo}
Here, $D\phi_s(x)$ denotes the derivative of the map $(s,x)\mapsto \phi(s,x)$
with respect to the variable $x$. We make some basic technical assumptions 
(continuity, Lipschitz character on the function $\phi$
with respect to the second variable)
to make sure that the flow is well defined.
We refer the reader to Section~\ref{se12} for more information,
and in particular to Theorem~\ref{theo-nondp2} below, where 
the definition of the contraction rate can be found.

The idea
of the proof is to construct a special flow-invariant set,
appealing to a generalization due to Martin~\cite{MR0318991}
 of theorems of Bony~\cite{bony} and Brezis~\cite{MR0257511} on
the geometric characterization of flow invariance. 
Formula~\eqref{e-explicit} should be compared with
results of Nussbaum, who studied the more general
question of computing the contraction rate
of a non-necessarily order-preserving flow in Thompson metric~\cite{MR1269677},
and obtained an explicit formula 
in the specical case of the standard positive cone.
However, this formula, valid for non order-preserving flows,
appears to have no natural generalization to abstract
cones (although some reasonably explicit
conditions can be given in the special case of symmetric cones, we
 leave this for a further work, see also Section~\ref{section-loss} below for
a special case). In addition, Nussbaum's
approach, which relies on the Finsler structure of the Thompson metric,
is widely applicable in its spirit but leads to different technical assumptions,
including geodesic convexity. 
See~\S\ref{sec-nussbaum} for a detailed comparison.



Then, we show, in Section~\ref{subsec-std}, that the contraction results of 
Liverani and Wojtowski~\cite{liveraniw} and of Lawson and Lim~\cite{MR2338433}
concerning the standard Riccati equation~\eqref{e-standard}
with positive semi-definite matrices $\Sigma,Q$,
as well as new contraction results in the case 
when $\Sigma$ is not positive semi-definite, can be
recovered, or obtained, 
by an application of Formula~\eqref{e-explicit}.
This provides an alternative to the earlier approaches,
which relied on the theory of symplectic semigroups.
This will allow us to handle as well
situations in which the symplectic structure is missing,
as it is the case of the generalized Riccati differential equation.


Our second main result shows that 
the flow of the generalized Riccati differential
equation is a local contraction in Thompson metric.
\begin{theo}\label{th-2}
Assume that the coefficients of the generalized Riccati differential equation~\eqref{e-sriccnew} are constant,  and that the matrix $\left(\begin{smallmatrix}Q& L' \\L&R\end{smallmatrix}\right)$ is positive definite. Then, the flow of this equation is a strict contraction on the interior of the cone of positive definite matrices, and this contraction is uniform on any subset that is bounded from above in the Loewner order.
\end{theo}
This theorem follows from Theorem~\ref{localcontraction}
in Section~\ref{sertcgg}, where an explicit bound for the contraction rate
on an interval in the Loewner order is given. We shall also see in Section~\ref{sertcgg} that the flow of the generalized Riccati equation
is no longer a uniform contraction on interior of the cone,
which reveals a fundamental discrepancy with the case of the standard Riccati equation.
Then, motivated by earlier results of Chen, Moore, Rami, and Zhou (see~\cite{MR1819815} and \cite{MR1778371}) on the asymptotic behavior of the GRDE,
we identify (Theorem~\ref{gareNEW}) different assumptions under
which a trajectory of the GRDE converges exponentially to a stable
solution of the associated Generalized Algebraic Riccati Equation (GARE).
We also establish (Section~\ref{subsec-gen-disc}) analogous results
concerning the discrete time case. Then, we 
give a necessary and sufficient condition (Proposition~\ref{prop-discrete})
for the generalized discrete Riccati operator to be
a strict global contraction. 

Finally, in Section~\ref{section-loss}, we establish the following
negative result, which shows that the 
Thompson metric is essentially the only invariant
Finsler metric in which the flow of the GRDE is non-expansive
for all admissible values of the matrix data.
\begin{theo}\label{theo-final}
The flow of the generalized Riccati differential equation is non-expansive in the invariant Finsler metric arising
from a symmetric gauge function, regardless of the 
parameters $(A,B,C,D,L,Q,R)$, if and only if this symmetric gauge function is a scalar multiple of the sup-norm. 
\end{theo}
In particular, the flow of the GRDE is {\em not} non-expansive
in the invariant Riemannian metric, showing that Bougerol's
theorem on the contraction of the standard discrete Riccati
equation does not carry over to the GRDE.

\section{Preliminaries}
\subsection{Thompson's part metric}\label{subsec-def-th}
We first recall the definition and basic properties of Thompson's part metric.

Throughout the paper, 
$\cX$ is a real or complex Banach space with norm $|\cdot|$.
Let $\cX^*$ be the dual space of $\cX$. For any $x\in \cX$ and $q\in \cX^*$, denote by $\<q,x>$ the real part of $q(x)$:
\[
\<q,x>=\Real q(x) \enspace .
\]
Let $\C\subset \cX$ be a closed pointed convex cone, i.e., $\alpha \C\subset \C$ for $\alpha \in \R^{+}$, $\C+\C\subset \C$ and $\C\cap (-\C)=0$. The {\em dual cone} of $\C$ is defined by $$\C^*=\{z\in \cX^*:\braket{z, x}\geq0 \enspace \forall x\in \C\}\enspace.$$ 
We denote by $\C_0$ the interior of $\C$. We define the partial order $\leq $ 
induced by $\C$ on $\cX$ by
\[
x\leq y \Leftrightarrow y-x \in \C  
\]
so that
$$
x\leq y \Rightarrow \braket{z,x} \leq \braket{z,y},\enspace \forall z\in \C^*.
$$
We also define the relation $\ll$ by
\[
x\ll y \Leftrightarrow y-x \in \C_0
\enspace .
\]
For $x\leq y$ we define the order intervals:
$$
[x,y]:=\{z\in \cX|x\leq z\leq y\},\qquad 
(x,y):=\{z\in \cX|x\ll z\ll y\}.
$$
For $x\in \cX$ and $y\in \C_0$, following~\cite{nussbaum88}, we define
\begin{align}\label{a-eq4}
\begin{array}{l}
M(x/y):=\inf\{t\in \R:x\leq ty\}\\
m(x/y):=\sup\{t\in \R:x\geq ty\}
\end{array}
\end{align}
Observe
that since $y\in\mathcal{C}_0$, and since $\mathcal{C}$ is closed
and pointed, the two sets in~\eqref{a-eq4} are non-empty, closed, and bounded
from below and from above, respectively. In particular, $m$ and $M$ take
finite values.
\begin{defi} The \firstdef{Thompson part metric} between two elements $x$ and $y$ of $\C_0$ is
\begin{align}
d_T(x,y):=\log(\max \{M(x/y),M(y/x)\})
\label{e-def-th}
\enspace .
\end{align}
\end{defi}
It can be verified that $d_T(\cdot,\cdot)$ defines a metric on $\C_0$, namely for any $x,y,z\in \C_0$ we have
$$
d_T(x,y)\geq 0,\enspace d_T(x,y)=d_T(y,x),\enspace d_T(x,z)\leq d_T(x,y)+d_T(y,z), \enspace d_T(x,y)=0\Leftrightarrow x=y.
$$
A sufficient condition for $\C_0$ to be complete with respect to $d_T(\cdot,\cdot)$ is that $\C$ is a normal cone,  see~\cite{ThompsonMR0149237}.
We shall consider specially the case in which $\C$ is the cone of $n\times n$ positive semi-definite matrices. Then, it can be checked that, for all $A,B\in \C_0$, 
\[ M(A/B)= \max_{1\leq i\leq n} \log \lambda_i ,\enspace,
\]
where
$\lambda_1,\dots,\lambda_n$ are the eigenvalues of the matrix $B^{-1}A$
(the latter eigenvalues are real and positive)
so that the Thompson metric $d_T$ can be explicitly computed
from~\eqref{e-def-th}.
%

\subsection{Characterization of flow invariant sets}
We next recall some known results on the characterization of flow-invariant sets in terms of tangent cones, which will be used to characterize order-preserving non-expansive flows in Thompson's metric.


In the sequel, $J=[0,T)\subset \R$ is a possibly unbounded interval,
$\D\subset \cX$ is an open set,
and $\phi(t,x)$ is a function from $J\times \D$ to $\cX$. 
For $x\in \cX$ and $\cS\subset \cX$ we define the distance function: $$d(x,\cS)=\inf\{|x-y|:y\in \cS\}.$$


We study the following Cauchy problem:
\begin{align}\label{initialvaluepb}
\left\{
\begin{array}{l}
\dot x(t)=\phi(t,x(t)),
\\ x(s)=x_0.
\end{array}\right.
\end{align} 
By a solution of~(\ref{initialvaluepb}) on $[s,a)\subset J$ we mean a continuously differentiable function $t\mapsto x(t):[s,a)\rightarrow \D$ such that $x(s)=x_0$ and $\dot x(t)=\phi(t,x(t))$ for all $t\in[s,a)$.

Let $\cS$ be a closed subset of $\cX$. We say that the system $(\cS\cap \D,\phi)$ is \firstdef{flow-invariant} if every solution of (\ref{initialvaluepb}) leaves $\cS$ invariant, in the sense that for any $s\in J$ and $x_0\in \cS\cap \D$, the solution $x(t)$ must be in $\in \cS\cap \D$, for all $t\in [s,a)$.

Characterizations of flow invariant sets go back to the works of Bony~\cite{bony} and Brezis~\cite{MR0257511}. Several improvements, together with extensions to the infinite dimensional case can be found in~\cite{MR0303024}, \cite{MR0318991}, ~\cite{MR0367131} and~\cite{RedhefferandWalterMR0470401}. We shall actually need here an immediate consequence of a theorem of Martin~\cite{MR0318991}.

\begin{theo}[Theorem 1 of~\cite{MR0318991}]\label{MartinTheorem}
 Suppose that the following conditions hold:
\begin{itemize}
\item[(C1)] $\phi$ is a continuous function on $J\times \D$;
\item[(C2)] For every closed bounded set $K\subset \D$, there is a constant $L>0$ such that $$
|\phi(t,x)-\phi(t,y)|\leq L|x-y|,\quad \forall t\in J,x,y\in K \enspace;
$$
\item[(C3)] For all $t\in J$ and $x\in \cS\cap\D$, $$\lim_{h\downarrow 0} \frac{d(x+h\phi(t,x),\cS\cap\D)}{h}=0\enspace .$$
\item[(C4)] $\cS$ is convex.
\end{itemize}
Then the system $(\cS\cap\D,\phi)$ is flow-invariant.
\end{theo}

It is not difficult to prove that for $x\in \cS\cap \D$, $v\in \cX$ and sufficiently small $h>0$, $$d(x+hv,\cS)=d(x+hv,\cS\cap\D).$$ Thus, Condition (C3) is equivalent to:
\begin{itemize}
\item[(C5)] For all $t\in J$ and $x\in \cS\cap\D$, $$\lim_{h\downarrow 0} \frac{d(x+h\phi(t,x),\cS)}{h}=0.$$
\end{itemize}
Condition (C2) is a local Lipschitz condition for the function $\phi$, with respect to the second variable. Condition (C3) is a tangency condition (the vector field $\phi$ should not point outward the set $\cS\cap\D$).


\begin{defi}[Tangent cone~\cite{MR0367131}]\label{deftangent}
The {\em tangent cone} to a closed set $\cS\subset \cX$ at a point $x\in \cS$, written $T_{\cS}(x)$, is the set of vectors $v$ such that:
\begin{align}\label{deftan}
\liminf_{h\downarrow 0} \frac{d(x+hv,\cS)}{h}=0.
\end{align}
\end{defi}


\begin{rem}
When $\cS$ is a closed convex, we know that the limit in (\ref{deftan}) exists.
Thus, Condition (C5), equivalent to (C3) in Theorem~\ref{MartinTheorem}, can be replaced by:
\begin{itemize}
\item[(C6)] For all $t\in J$ and $x\in \cS\cap\D$,$$\phi(t,x)\in T_{\cS}(x).$$ 
\end{itemize}
\end{rem}




Besides, this definition coincides with the one in convex analysis, i.e.,
\begin{prop}[Proposition 5.5, Exercice 7.2~\cite{Clarke:1998:NAC:274798}]\label{convextangent}
Let $\cS$ be a closed convex set of $\cX$, then,
\begin{align*}
T_{\cS}(x)=\operatorname{cl}\{v:\exists \lambda >0~\text{with}~x+\lambda v\in \cS\},\enspace\forall x\in \cS.\end{align*}
\end{prop}

Now the flow-invariance can be checked by verifying if $\phi(t,x)$ lies in the tangent cone of $\cS$. To this end, we need to compute the tangent cone at each point of $\cS$. In some cases the tangent cone can be expressed in a simple way:
\begin{lemma}[Exercise 2.5.3~{\cite{Clarke:1998:NAC:274798}}]\label{tangentproduct}
 Let $\cS_1,\cS_2\subset \cX$ be closed subsets, $x=(x_1,x_2)\in \cS_1\times \cS_2$. Then
$$
T_{\cS_1\times \cS_2}(x)=T_{\cS_1}(x_1)\times T_{\cS_2}(x_2).
$$
\end{lemma}
We shall consider specially $\cS=\C$. Then, using Proposition~\ref{convextangent} and the Hahn-Banach separation theorem, one can show that
\begin{align}\label{TCX}
T_{\C}(x)=\{v| \<q,v>\geq 0 \enspace \text{if} \enspace q\in \C^* \enspace \text{and}\enspace \<q,x>=0\},\enspace x\in \C.
\end{align}

\section{The contraction rate in Thompson metric of order-preserving flows}\label{se12}


\subsection{Preliminary results}
From now until the end, the function $\phi(t,x)$ is assumed to be continuous on $J\times \D$ and Fr\'echet differentiable to $x$. The derivative of $\phi$ with respect to the second variable at point $(t,x)$ is denoted by $D\phi_t(x)$. We also assume that the derivative is bounded on any closed bounded set, i.e., for any bounded set $K\subset \D$, there is a constant $L$ such that:
$$
|D\phi_t(x)|\leq L,\quad\forall t\in J, x\in K.
$$
Therefore Condition (C1) and Condition (C2) are both satisfied.
The existence and uniqueness of the solution of~(\ref{initialvaluepb}) follow from the Cauchy-Lipschitz Theorem. We then define the flow $M_{\cdot}^{\cdot}(\cdot)$ associated to the system by:
$$
M_s^t(x_0)=x(t), \enspace s\in J,\enspace t\in [s,a)
$$
where $x(t):t\in [s,a)$ is the maximal solution of~(\ref{initialvaluepb}).
(Note that in general the flow is defined only on a subset of $J\times J\times \D$.)
For each open subset $\cU\subset \D$ and initial value $x_0\in \cU$ we define $t_\cU(s,x_0)$ as the first time when the trajectory leaves $\cU$, i.e.,
$$
t_\cU(s,x_0)=\sup\{b\in (s,a)| M_s^t(x_0)\in \cU, \enspace \forall t\in [s,b)\}.
$$
When $\phi$ is independent of time $t$, we denote simply
$$M_{t}(x_0):=M_0^t(x_0),\quad t_{\cU}(x_0):=t_{\cU}(0,x_0),\enspace \forall x_0\in \D, \cU\subset \D.$$
By uniqueness of the solution, the flow has the group property:
$$
M_{s}^{t}(x_0)=M_{s}^{t_1}(M_{t_1}^{t}(x_0)),\quad\forall \enspace 0\leq  s \leq t_1\leq t<t_{\D}(s,x_0)
\enspace .
$$

\begin{defi}[Order-preserving flow]\label{monot}
Let $\cU$ be an open subset of $\D$. The flow $M_{\cdot}^{\cdot}(\cdot)$ is said to be \textit{order-preserving on $\cU$} if for all $x_1,x_2 \in \cU$ such that $x_1\geq x_2$,
$$ M_s^t(x_1)\geq M_s^t(x_2),\enspace \forall \enspace 0\leq s \leq t< t_{\cU}(s,x_1)\wedge t_{\cU}(s,x_2).$$
\end{defi}

\begin{defi}[Non-expansiveness and contraction]\label{defnonexp}
Suppose that $\C_0\subset \D$. The flow $M_{\cdot}^{\cdot}(\cdot)$ is said to be {\em contractive} on $\C_0$ 
with rate $\alpha>0$ in Thompson metric if for all $x_1,x_2\in \C_0$,
$$
d_T(M_s^t(x_1),M_s^t(x_2))\leq e^{-\alpha (t-s)} d_T(x_1,x_2),\enspace \forall \enspace 0\leq s \leq t< t_{\C_0}(s,x_1)\wedge t_{\C_0}(s,x_2)
$$
If the latter inequality holds with $\alpha=0$, the flow is said to be 
{\em non-expansive}.
\end{defi}

In the following, our primary goal is to characterize non-expansive order-preserving flows in Thompson part metric.
 We shall need the following proposition, which provides a characterization of monotonicity in terms of the function $\phi$.  
The equivalence of the first two assertions was proved in~\cite{RedhefferandWalterMR0470401}.

\begin{prop}[Compare with Theorem~3 in~\cite{RedhefferandWalterMR0470401}]\label{croissance}
 Let $\cU$ be an open subset of $\D$. The following conditions are equivalent:
\begin{itemize}
\item[(a)] The flow $M_{\cdot}^{\cdot}(\cdot)$ is order-preserving on $\cU$.
\item[(b)] For all $s\in J$ and $x_1,x_2\in \cU$ such that $x_1\geq x_2$, $\phi(s,x_1)-\phi(s,x_2)\in T_{\C}(x_1-x_2)$.
\end{itemize}
 If $\cU$ is convex, then the above conditions are equivalent to:
\begin{itemize}
\item[(c)]For all $s\in J$, $x\in \cU$ and $v\in \C$,
\begin{align}\label{condpro}
\braket{q, D\phi_s(x) v}\geq 0,\enspace \forall q\in \{q\in  \C^*:\braket{q, v}=0\}.
\end{align}
\end{itemize}

\end{prop}
\begin{proof}


We only need to prove the equivalence between (b) and (c), since
the equivalence between (a) and (b) follows from~\cite{RedhefferandWalterMR0470401}.
In view of~\eqref{TCX}, Condition (b) is equivalent to
the following:

for all $s\in J$ and $x_1,x_2\in \cU$ such that $x_1\geq x_2$,
 $$\<q,\phi(s,x_1)-\phi(s,x_2)> \geq 0,\enspace \forall q\in\{q\in \C^*:\<q,x_1-x_2>=0\}.$$
Now suppose that (b) is true. Then for any $s\in J$, $x\in \cU$ and any $v\in \C$, there is $\delta >0$ such that for any $0\leq \epsilon\leq \delta$
$$
\<q,\phi(s,x+\epsilon v)-\phi(s,x)>\geq 0,\enspace \forall q\in\{q\in \C^*:\<q,v>=0\}.
$$
Since $\phi$ is differentiable at point $x$, dividing by $\epsilon$ the latter
inequality, and letting $\epsilon$ tend to $0$, we get
$$
\<q,D\phi_s(x) v>\geq 0,\enspace \forall q\in\{q\in \C^*:\<q,v>=0\}.
$$
Next suppose that Condition (c) holds.
Fix any $s\in J$ and $x_1,x_2\in \cU$ such that $x_1\geq x_2$. Fix any $q\in \C^*$ such that $\<q,x_1-x_2>=0$. Define the function $g:[0,1]\rightarrow \R$ by:
$$
g(\lambda)=\<q,\phi(s,\lambda x_1+(1-\lambda)x_2)-\phi(s,x_2)>.
$$
Then we have $g(0)=0$ and in view of convexity of $\cU$ and \eqref{condpro},
$$g'(\lambda)=\braket{q, D\phi_s(\lambda x_1+(1-\lambda) x_2) (x_1-x_2))}\geq 0,\enspace \forall 0\leq \lambda \leq 1.$$
 A standard argument establishes that:
\begin{align}\label{croz}
g(1)=\braket{q, \phi(s,x_1)-\phi(s,x_2)}\geq 0.
\end{align}
Since $s,x_1,x_2$ and $q$ are arbitrary, we deduce Condition (b).
\end{proof}

\subsection{Characterization of the contraction rate in terms of flow invariant sets}
The following is a key technical result in the characterization of the contraction rate of the flow. 

\begin{prop}\label{nondp2}
Let $\cU\subset \D$ be an open set such that $\lambda\cU\subset \cU$ for all $\lambda\in (0,1]$.
 If the flow $M_{\cdot}^{\cdot}(\cdot)$ is order-preserving on $\cU$, then the following conditions are equivalent:
\begin{itemize}
 \item [(a)]For all $x\in \cU$ and $\lambda\geq 1$ such that $\lambda x\in \cU$,
$$
M_s^t(\lambda x)\leq \lambda^{e^{-\alpha (t-s)}}M_s^t(x),\enspace  0\leq s\leq t< t_{\cU}(s,x)\wedge t_{\cU}(s,\lambda x).
$$

\item [(b)] For all $s\in J$ and $x\in \cU$,
$$D\phi_s(x)x-\phi(s,x)\leq -\alpha x. 
$$
\item [(c)] For all $x, y\in \cU$ and $\lambda \geq 1$ such that $y\leq \lambda x$,
$$
M_s^t(y)\leq \lambda^{e^{-\alpha (t-s)}}M_s^t(x), 0\leq s\leq t< t_{\cU}(s,x)\wedge t_{\cU}(s,y).
$$
\end{itemize}
\end{prop}
\begin{proof}
Suppose Condition (a) holds. Let $x$ be any point in $\cU$. Fix any $\lambda > 1$ such that $\lambda x\in \cU$, we must have:
 \begin{align}\label{eere}
M_s^t(\lambda x)\leq \lambda^{e^{-\alpha (t-s)}}M_s^t(x),\qquad  0\leq t< t_{\cU}(s,x)\wedge t_{\cU}(s,\lambda x).
\end{align}
where it must be the case that $t_{\cU}(s,x)\wedge t_{\cU}(s,\lambda x)>0$. 
Since the terms on both sides of~\eqref{eere} coincide when $t=s$,
taking the derivative of each of these terms at $t=s$, we obtain
\begin{align}\label{ep}
 \phi(s,\lambda x)\leq \lambda \phi(s,x)-\alpha (\lambda \ln \lambda) x
\end{align}
Since this inequality holds for all $\lambda\geq 1$ such that $\lambda x\in \cU$, with equality for $\lambda=1$, the derivation of the two sides of the above inequality at $\lambda=1$ leads to:
$$
D\phi_s(x) x-\phi(s,x)\leq -\alpha x 
$$
for all $x\in \cU$. Condition (b) is deduced.

Now suppose that Condition (b) is true. We shall derive Condition~(c)
by constructing an invariant set. 
Denote:
$$
\tilde \cX:=\cX\times  \cX\times \R,$$
$$
\tilde \D:=\cU\times \cU\times \R^+\backslash\{0\},
$$$$
 \cS:=\{(x_1,x_2,\lambda)\in \tilde \cX: x_2\leq \lambda x_1,\lambda\geq 1\}.
$$
Define the differential equation on $\tilde\D$:
\begin{align}\label{equ1old}
\left(
\begin{array}{l}
\dot x_1\\ \dot x_2 \\ \dot \lambda
\end{array}\right)=\Phi(t,x_1,x_2,\lambda):=\left(
\begin{array}{l}
 \phi(t,x_1)\\ \phi(t,x_2)\\ -\alpha \lambda \ln \lambda
\end{array}
\right)
\end{align}
It is not difficult to see that Condition (c) is equivalent to the flow-invariance of the system $(\cS\cap \tilde\D,\Phi)$. 
It would be natural to show directly the latter flow-invariance by appealing to Theorem~\ref{MartinTheorem}, but
the set $\cS$ is not convex, making it harder to check the assumptions of this theorem. Therefore, we make a change of variable to replace
$\cS$ by a convex set.

Define the smooth function $F:\tilde \cX \rightarrow \tilde\cX$ by:
$$
F(x_1,x_2, \lambda)=(x_1,\lambda x_1-x_2,\lambda-1),\enspace\forall (x_1,x_2,\lambda)\in \tilde \cX.
$$
Denote $$\cS'=\cX\times \C\times \R^+.$$ By Lemma~\ref{tangentproduct}, for $(y_1,y_2,\kappa)\in \cS'$,
$$
T_{\cS'}(y_1,y_2,\kappa)=\cX\times T_{\C}(y_2)\times T_{\R^+}(\kappa).
$$
Observe that $\cS=\{x\in \tilde \cX| F(x)\in \cS'\}$ and that $F$ has a smooth inverse $G:\tilde \cX\rightarrow \tilde \cX$ given by:
$$G(y_1,y_2,\kappa)=(y_1,(\kappa+1)y_1-y_2,\kappa+1).$$ 
Therefore $F(\tilde\D)=G^{-1}(\tilde\D)$ is an open set.
Let $(y_1,y_2,\kappa)=F(x_1,x_2,\lambda)$ and consider the system:
\[
(\dot y_1, \dot y_2, \dot \kappa)'
=\Psi
(t,y_1,y_2,\kappa)
\]
where
\[
\Psi
(t,y_1,y_2,\kappa)
:=\left(
\begin{array}{l}
 \phi(t,y_1)\\ -\alpha (\kappa+1) \ln (\kappa+1) y_1+(\kappa+1) \phi(t,y_1)-\phi(t,(\kappa+1) y_1-y_2)\\-\alpha (\kappa+1) \ln (\kappa+1)
\end{array}
\right)
\]
One can verify that the invariance of the system $(\cS\cap \tilde\D,\Phi)$ is equivalent to the invariance of the system $(\cS'\cap  F(\tilde\D),\Psi)$.

Now the function $\Psi:J\times F(\tilde \D)\rightarrow \tilde \cX$ defined as above is continuous and differentiable to the second variable with bounded derivative on bounded set. Besides $\cS'$ is convex.
By applying Theorem~\ref{MartinTheorem}, the system $(\cS'\cap F(\tilde\D),\Psi)$ is flow-invariant if the following condition is satisfied:
\begin{align}\label{tangencycondition}
\Psi(s,y_1,y_2,\kappa)\in T_{\cS'}(y_1,y_2,\kappa),\enspace\forall s\in J, (y_1,y_2,\kappa)\in \cS'\cap F(\tilde\D).
\end{align}
That is, for any $(y_1,y_2,\kappa) \in \cS'\cap F(\tilde\D)$ and $s\in J$,
$$\left\{\begin{array}{l}
\phi(s,y_1)\in \cX \\
-\alpha (\kappa+1) \ln (\kappa+1) y_1+(\kappa+1) \phi(s,y_1)-\phi(s,(\kappa+1) y_1-y_2)\in T_{\C}(y_2)\\
-\alpha (\kappa+1) \ln (\kappa+1)\in T_{\R^+}(\kappa)
\end{array}\right.
$$
It suffices to check the second condition because the others hold trivially. By applying the bijection $F$, this condition becomes: for any $s\in J$ and $(x_1,x_2,\lambda)\in \cS\cap \tilde\D$,
$$
-\alpha \lambda \ln \lambda x_1+\lambda \phi(s,x_1)-\phi(s,x_2)\in T_{\C}(\lambda x_1-x_2).
$$
Let any $s\in J$, $x_1,x_2\in \cU$ and 
$\lambda\geq 1$ such that $x_2\leq \lambda x_1$. Let any $q\in \C^*$ such that $\<q,\lambda x_1-x_2>=0$. By~(\ref{TCX}) we only need to prove:
\begin{align}\label{ester}
\<q,-\alpha \lambda \ln \lambda x_1+\lambda \phi(s,x_1)-\phi(s,x_2)>\geq 0.
\end{align}
By the assumptions, we know that $\lambda^{-1}x_2\in \cU$. Then, it suffices to prove: for any $x_1,x_2\in \cU$ such that $x_1\geq x_2$, let $q\in \C^*$ such that $\<q, x_1-x_2>=0$, then for any $\lambda\geq 1$ such that $\lambda x_2\in \cU$ we have:
$$
\<q,-\alpha \lambda \ln \lambda x_1+\lambda \phi(s,x_1)-\phi(s,\lambda x_2)>\geq 0.
$$
Define the function $f: [1,\lambda]\rightarrow \R$ by:
$$
f(\tau)=\braket{q,-\alpha \ln \tau x_1+\phi(s,x_1)-\tau^{-1}\phi(s,\tau x_2)}
$$
Notice that the function $f$ is well defined on $[1,\lambda]$.
By hypothesis of monotonicity and Proposition~\ref{croissance},
$$
f(1)=\braket{q,\phi(s,x_1)-\phi(s,x_2)}\geq 0.
$$
Differentiating $f$ gives,
for all $\tau \in[1,\lambda]$,
\begin{align*}
f'(\tau)&=\braket{q,-\tau^{-1}\alpha x_1+\tau^{-2}\phi(s, \tau x_2)-\tau^{-1}D\phi_s(\tau x_2) x_2} \\ 
&\geq  \braket{q,-\tau^{-1}\alpha x_1+\tau^{-1}\alpha x_2}  
\qquad \mbox{(by Condition (b))}\\
& = 0  \enspace .
\end{align*}
A standard argument establishes that
$f(\lambda)\geq 0$, and so~\eqref{ester} is proved, whence the flow-invariance of $(\cS\cap \tilde\D,\Phi)$, which is exactly Condition (c).  Finally, Condition (a) follows from Condition (c) by considering $y=\lambda x$.
\end{proof}

We next state the main results. Recall that $J=[0,T)\subset \R$.
\begin{theo}[Contraction rate]\label{theo-nondp2}
Assume that $\phi$ is defined on $J\times \cU$ where $\cU \subset\C_0$ is an open set in the interior of the cone
such that $\lambda\cU\subset \cU$ for all $\lambda\in (0,1]$.
If the flow $M_{\cdot}^{\cdot}(\cdot)$ is order-preserving on $\cU$, then
the best constant $\alpha$ such that
\begin{align}\label{contheo}
d_T(M_s^t(x_1),M_s^t(x_2))\leq e^{-\alpha (t-s)}d_T(x_1,x_2), \enspace 0\leq s\leq t< t_{\cU}(s,x_1)\wedge t_{\cU}(s,x_2)
\end{align}
holds for all $x_1,x_2\in \cU$ is given by 
\begin{align}\label{formulaal}
\alpha:=-\sup_{s\in J,\; x\in \cU}  M\big((D\phi_s(x)x-\phi(s,x))/x\big) \enspace .
\end{align}
\end{theo}
\begin{proof}
If~\eqref{contheo} holds for all $x_1, x_2\in \cU$, then Condition (a) in Proposition~\eqref{nondp2} holds. It follows that the constant $\alpha$ must satisfy
\begin{align}\label{inealpha}
D\phi_s(x)x-\phi(s,x)\leq -\alpha x, \enspace \forall s\in J, x\in \cU.
\end{align}
Now conversely if~\eqref{inealpha} holds. Then Condition (c) in
Proposition~\ref{nondp2} holds. For any $x_1,x_2\in \cU$,
let $\lambda=e^{d_T(x_1,x_2)}$, then
$$
M_s^t(x_1)\leq \lambda^{e^{-\alpha(t-s)}}M_s^t(x_2), 0\leq s\leq t<t_{\cU}(s,x_2)\wedge t_{\cU}(s,x_1) \enspace .
$$ 
The same is true if we exchange the roles of $x_1$ and $x_2$,
and so,
~\eqref{contheo} holds for all $x_1,x_2\in \cU$. Consequently the best constant $\alpha$ such that~\eqref{contheo} holds for all $x_1,x_2\in \cU$ must be the greatest constant $\alpha$ such that~\eqref{inealpha} holds, which 
is precisely~\eqref{formulaal}.
\end{proof}

Now we get a direct corollary.
\begin{theo}\label{nonexp}
Suppose that $\phi$ is defined on $J\times \C_0$. Let $\alpha \in \R$.
If the flow is order-preserving on $\C_0$, then the following are equivalent:
\begin{enumerate}
\item [(a)]For all $x_1,x_2\in \C_0$:
$$
d_T(M_s^t(x_1),M_s^t(x_2))\leq e^{-\alpha (t-s)}d_T(x_1,x_2), \enspace 0\leq s\leq t< t_{\C_0}(s,x_1)\wedge t_{\C_0}(s,x_2).
$$
\item [(b)]For all $s\in J$ and $x\in \C_0$,$$\label{p2}
D\phi_s(x)x-\phi(s,x)\leq -\alpha x. 
$$
\end{enumerate}
If any of these conditions holds, then the flow leaves $\C_0$ invariant, i.e., for any $s\in J$ and $x\in \C_0$, $t_{\C_0}(s,x)=T$.
\end{theo}
\begin{proof}
The equivalence between (a) and (b) follows from Theorem~\ref{theo-nondp2}.
Now suppose that Condition (b) holds. Let any $s\in J$ and $x_1,x_2\in \C_0$. Let $t_1=t_{\C_0}(s,x_1)$ and $t_2=t_{\C_0}(s,x_2)$. 
Suppose that $t_1<t_2$. Then it must be the case that $t_1<+\infty$. Thus the set 
$\{M_s^r(x_2):r\in [s,t_1]\}$
is compact and included in $\C_0$. Denote
$$
K=\max \{d_T(M_s^r(x_2),M_s^{t_1}(x_2)) | \enspace r\in [s,t_1]\}<+\infty
$$
and $K_0=K+\max \{e^{-\alpha (t_1-s)},1\}d_T(x_1,x_2) $.
Note that there exists $s<\bar r<t_1$ such that 
$$
d_T(M_s^{\bar r}(x_1),M_s^{t_1}(x_2))>K_0,
$$
otherwise $t_{\C_0}(s,x_1)>t_1$. 
But for any $s<r<t_1$,
$$
\begin{array}{ll}
d_T(M_s^r(x_1),M_s^{t_1}(x_2))&\leq d_T(M_s^r(x_1),M_s^r(x_2))+d_T(M_s^r(x_2),M_s^{t_1}(x_2)) \\&\leq e^{-\alpha (r-s)}d_T(x_1,x_2)+d_T(M_s^r(x_2),M_s^{t_1}(x_2))\\&\leq K_0.
\end{array}
$$
The contradiction implies that $t_1< t_2$ is impossible. We then showed that there exists $\bar T\in (0,+\infty]$ such that for any $s\in J$ and $x\in \C_0$, $t_{\C_0}(s,x)=\bar T$. 
From the group property of the flow action, we deduce that $\bar T=T$.
\end{proof}

In the sequel we suppose that the dynamics $\phi$ is independent of time and study the convergence of an orbit of the flow to a 
fixed point in the interior of the cone. Let $\bar x\in \C_0$ be such that $\phi(\bar x)=0$. Let $\mu>1$. Denote by $\cU$ the open interval $(\mu^ {-1}\bar x,\mu \bar x)$. We look for the best constant $\alpha \in \R$ such that:
\begin{align}\label{agryfg}
d_T(M_t(x), \bar x)\leq e^{-\alpha t} d_T(x,\bar x),\enspace\forall x \in \cU, 0\leq t< t_{\cU}(x).
\end{align}
\begin{theo}[Convergence rate]\label{tzsdfg}
 We assume that $\phi$ is independent of time, defined on $\C_0$ and such that the flow is order-preserving on $\C_0$. Let $\bar x\in \C_0$ be a zero point of $\phi$. Then the best constant $\alpha$ such that~(\ref{agryfg}) holds 
is given by
\begin{align}\label{rtsdgfnew}
\alpha = \inf_{\mu^{-1}< \lambda < \mu} m\big((-(\lambda \ln \lambda)^{-1} \phi(\lambda \bar{x}) )/ \bar{x}\big) \enspace .
\end{align}
Moreover, if the latter $\alpha$ is non-negative, then for all $x\in [\mu^ {-1}\bar x, \mu \bar x]$,
\begin{align}\label{agryfger}
d_T(M_t(x), \bar x)\leq e^{-\alpha t} d_T(x,\bar x),\enspace\forall t\geq 0.
\end{align}
\end{theo}
\begin{proof}
Suppose that $\alpha$ satisfies~\eqref{agryfg}.
 Let any $\lambda\in (1,\mu)$. Then $\lambda \bar x \in \cU$ and
$$
M_t(\lambda \bar x)\leq \lambda^ {e^ {-\alpha t}}\bar x ,\enspace  0\leq t<t_{\cU}(\lambda \bar x).
$$
Since $t_{\cU}(\lambda \bar x)>0$ and both sides of the former inequality
coincide when $t=0$, we get the inequality for the derivative at $t=0$:
\begin{align}\label{atuefg}
\phi(\lambda \bar x)\leq -\alpha \lambda (\ln \lambda) \bar x \enspace,
\end{align}
and so 
\begin{align}
 \alpha \bar x\leq -(\lambda \ln \lambda)^ {-1}\phi(\lambda \bar x), \enspace \forall 1<\lambda <\mu.
\end{align}
Similarly, for $\lambda \in (\mu^ {-1},1)$,
$$
\lambda^ {e^ {-\alpha t}} \bar x \leq M_t(\lambda \bar x),\enspace 0\leq t<t_\cU(\lambda \bar x),
$$ 
thus
\begin{align}\label{aerdfg}
- \alpha \lambda  \ln \lambda \bar x\leq \phi(\lambda \bar x)
\end{align}
leading to
\begin{align}
 \alpha \bar x\leq -(\lambda \ln \lambda)^ {-1}\phi(\lambda \bar x), \enspace \forall \mu^{-1}<\lambda <1.
\end{align}
It follows that $\alpha$ is bounded above by the expression in~\eqref{rtsdgfnew}.
To prove that conversely, ~\eqref{agryfg} holds when $\alpha$ is given
by~\eqref{rtsdgfnew}, we use an invariance argument as in the proof of Proposition~\ref{nondp2}. 
Denote:
$$
\tilde \cX:=  \cX\times \R,$$
$$
\sD:=\cU\times (1,\mu),
$$
$$
\cS_1:=\{(x,\lambda)\in \tilde \cX:  x\leq \lambda \bar x\},
$$
$$
\cS_2:=\{(x,\lambda)\in \tilde \cX:\bar x \leq \lambda x\},
$$
and define the differential equation:
\begin{align}\label{equ1}
\left(
\begin{array}{l}
\dot x \\ \dot \lambda
\end{array}\right)=\Phi(x,\lambda):=\left(
\begin{array}{l}\phi(x)\\ -\alpha \lambda \ln \lambda
\end{array}
\right).
\end{align}
Then ~(\ref{agryfg}) holds if $(\cS_1\cap \sD,\Phi)$ and $(\cS_2\cap \sD,\Phi)$ are invariant systems. Given the convexity of $\cS_1$, we can directly apply Theorem~\ref{MartinTheorem} to prove the invariance of the system $(\cS_1\cap \sD,\Phi)$. The tangent cone of $\cS_1$ at point $(x,\lambda)\in \cS_1$ is given by:
$$
T_{\cS_1}(x,\lambda)=\{(z,\eta): \<q,\eta \bar x-z>\geq 0, \forall  q\in \C^ {*},\<q,\lambda \bar x-x>=0\}.
$$
For any $q\in \C^ {*}$ such that $\<q,\lambda \bar x-x>=0$, by the order-preserving assumption and Proposition~\ref{croissance},
$$\<q,\phi(\lambda \bar x)>\geq \<q,\phi(x)>.$$
Now, using the expression of $\alpha$ in~\eqref{rtsdgfnew},
$$
\<q,-\alpha \lambda \ln \lambda \bar x-\phi(x)>\geq \<q,-\alpha \lambda \ln \lambda \bar x-\phi(\lambda \bar x)>\geq 0.
$$
This shows that $$\Phi(x,\lambda)\in T_{\cS_1}(x,\lambda),\enspace \forall (x,\lambda)\in \sD\cap \cS_1,$$whence the invariance of $(\cS_1\cap \sD,\Phi)$.
For the invariance of system $(\cS_2\cap \sD,\Phi)$, we define a bijection on $\sD$:
$$
F(x,\lambda)=(\lambda x-\bar x,\lambda)
$$
whose inverse is:
$$
G(y,\kappa)=( \kappa^{-1}(\bar x +y),\kappa).
$$
If $(x(\cdot),\lambda(\cdot))\in \sD$ follows the dynamics of~(\ref{equ1}), then $(y(\cdot),\kappa(\cdot))=F(x(\cdot),\lambda(\cdot))$ is the solution of the following differential equation:
\begin{align}\label{equ1-quat}
\left(
\begin{array}{l}
\dot y \\ \dot \kappa
\end{array}\right)=\Psi(y,\kappa)=\left(
\begin{array}{l}-\alpha\ln\kappa(\bar x+y)+\kappa \phi(\kappa^ {-1}(\bar x+y))\\ -\alpha \kappa \ln \kappa
\end{array}
\right).
\end{align}
Thus the invariance of system $(F(\sD)\cap F(\cS_2),\Psi)$ implies the invariance of system $( \sD\cap \cS_2,\Phi)$. Note that $F(\cS_2)=\C\times \R$. Therefore by Theorem~\ref{MartinTheorem} the system $( F(\sD)\cap F(\cS_2),\Psi )$ is invariant if
$$
\Psi(y,\kappa)\in T_{F(\cS_2)}(y,\kappa),\enspace\forall (y,\kappa)\in F(\sD)\cap F(\cS_2).
$$
The tangent cone of $F(\cS_2)$ at point $(y,\kappa)\in F(\cS_2)$ is given by:
$$
T_{F(\cS_2)}(y,\kappa)=\{z: \<q,z>\geq 0,\forall q\in \C^ {*}, \<q,y>=0\}\times \R.
$$
Again by the order-preserving assumption, for any  $q\in \C^ {*}$ such that $\<q,y>=0$,
$$
\<q, \phi(\kappa^ {-1}(\bar x+y))>\geq \<q,\phi(\kappa^ {-1}(\bar x))> 
$$
Using again the expression of $\alpha$ in~\eqref{rtsdgfnew},
$$
\<q,\kappa\phi(\kappa^{-1}(\bar x))>\geq \<q,(\alpha \ln \kappa)\bar x>
$$
because $\kappa \in (1,\mu).$ Therefore
$$
\<q,-(\alpha \ln \kappa )(\bar x+y)+\kappa \phi(\kappa^ {-1}(\bar x+y))>\geq 0,
$$
which implies
$$\Psi(y,\kappa)\in T_{F(\cS_2)}(y,\kappa),\enspace \forall (y,\kappa)\in F(\sD)\cap F(\cS_2),$$whence the invariance of $(F(\cS_2)\cap F(\sD),\Psi)$ and that of $(\cS_2\cap \sD,\Phi)$. 

Finally, if $\alpha\geq 0$, then the set $\cU$ is invariant (by~\eqref{agryfg}). Thus $t_{\cU}(x)=+\infty$ for all $x\in \cU$. Since the closure $[\mu^ {-1} \bar x, \mu \bar x]$ of $\cU$ is in the interior of the cone, we conclude that the relation~\eqref{agryfg} holds as well for $x\in [\mu^ {-1} \bar x, \mu \bar x]$.
\end{proof}

\subsection{The discrete time case}
For completeness, we give in this section the results analogous to Proposition~\ref{croissance} and~Theorem~\ref{theo-nondp2} for discrete operators, which
are of a simpler character. In this section we consider a differentiable map $F:\C_0\rightarrow \C_0$.
The first proposition characterizes order-preserving maps, its elementary
proof is left to the reader. 
\begin{prop}\label{pazeff}
 Let $\cU\subset \C_0$ be any open convex set. Then $F$ is order-preserving on $\cU$ if and only if
$$
DF(P)\cdot Z\geq 0,\enspace\forall P\in \cU, Z\in \C
$$ 
\end{prop}
Let $\sG\subset \C_0$. The Lipschitz constant of $F$ on $\sG$, denoted by $\mylip{F;\sG}$, is defined as:
\begin{align}\label{ruapfg-nocited}
\mylip{F;\sG}:=\displaystyle\sup_{P_1,P_2\in \sG} \frac{d_T(F(P_1),F(P_2))}{d_T(P_1,P_2)}.
\end{align}
\begin{prop}\label{ptthsd}
 Let $\sG\subset \C_0$ be a set such that $t\sG\subset \sG$ for any $t\geq 1$. If $F$ is order-preserving on $\sG$, then $$\mylip{F;\sG}=\inf\{\alpha:DF(P)\cdot P\leq \alpha F(P),\enspace \forall P\in \sG\}.$$
\end{prop}
\begin{proof}
It suffices to prove the equivalence between the following two conditions:
\begin{itemize}
 \item [(a)] $d_T(F(P_1),F(P_2))\leq \alpha d_T(P_1,P_2),\enspace \forall P_1,P_2\in \sG$
\item [(b)] $DF(P)\cdot P\leq \alpha F(P),\enspace \forall P\in \sG$
\end{itemize}
As was pointed out in Remark 1.9~\cite{MR1269677}, if $F$ is order-preserving, then Condition (a) is true if and only if:
$$
\lambda^ {-\alpha}F(\lambda P)\leq F(P),\forall P\in \sG,\lambda \geq 1.
$$ 
Condition (b) is a necessary condition (differentiate the above inequality at $\lambda=1$). For the sufficiency, note that the derivative of the left-hand side is
$$
\lambda^ {-\alpha-1}(DF(\lambda P)\cdot(\lambda P)-\alpha F(\lambda P))
$$
which is always negative semi-definite given that Condition (b) is true.
\end{proof}
\begin{rem}\label{rem-er}
Nussbaum treated the discrete case in~\cite{MR1269677}, as an intermediate step before considering differential equations. Corollary 1.3 there shows that for any open subset $\sG\subset \C_0$ such that for all $u,v\in \sG$ there exists a piecewise $\mathscr{C}^1$ minimal geodesic contained in $\sG$ (geodesic convexity assumption), the Lipschitz constant of the map $F$ on $\sG$ satisfies :
\begin{align}\label{e-twoside}
\mylip{F;\sG}=\inf\{\alpha: -\alpha F(P)\leq DF(P)\cdot Z\leq \alpha F(P), \enspace \forall P\in \sG, -P\leq Z\leq P \}
\end{align}
Thus, when the map $F$ is order-preserving, a variant of Proposition~\ref{ptthsd},
in which the domain $\sG$ satisfies the previous geodesic convexity assumption
can be easily obtained as a corollary of this result.
%
\end{rem}

\section{First applications and illustrations}\label{sec-first}

In this section, we show that several known contraction results, which were originally obtained in~\cite{liveraniw} and \cite{MR2338433} by means of symplectic semigroups, as well
as new ones concerning the standard Riccati equation with indefinite coefficients, can be obtained readily from Theorem~\ref{theo-nondp2}. The extension of these results to the generalized Riccati equation will be dealt with in Section~\ref{sec-sto}.
\subsection{Contraction rate of order-preserving flows on the standard positive cone}\label{sec-nussbaum}
Let us consider the standard cone $\C:=\{x\in \R^n:x_i\geq 0,1\leq i\leq n\}$ in $\cX=\R^n$ and an order-preserving flow $M_{\cdot}^{\cdot}(\cdot)$ associated to a differentiable function $\phi:J\times \C\rightarrow \R^n$. For a subset $\cU\subset \C_0$, define the best contraction rate on $\cU$ to be the greatest value of $\alpha$ satisfying: 
\begin{align}\label{leastcontractionrate}
 d_T(M_s^t(x_1),M_s^t(x_2))\leq e^{-\alpha (t-s)}d_T(x_1,x_2),\forall x_1,x_2\in \cU, 0\leq s\leq t < t_\cU(s,x_1)\wedge t_\cU(s,x_2)
\end{align}

A direct application of Theorem~\ref{theo-nondp2} is the following:
\begin{coro}[Compare with~{\cite[Th.~3.10]{MR1269677}}]\label{coroNus}
 Let $\cU\subset \C_0$ be an open set satisfying $\lambda \cU\subset \cU$ for all $\lambda \in (0,1]$.
For $s\in J$ and $x\in \cU$ define $g_i(s,x)$ by:
\begin{align}
g_i(s,x)=-x_i^{-1}[\sum_{j=1}^n \frac{\partial \phi_i}{\partial x_j}(s,x)x_j-\phi_i(s,x)]\label{e-simpl}
\end{align}
then the best contraction rate on $\cU$ defined in~\eqref{leastcontractionrate} is given by:
\begin{align}\label{atdgh}
\alpha=\displaystyle\inf \{g_i(s,x):1\leq i\leq n, x\in \cU,s\in J\}
\end{align}
\end{coro}
This should be compared with a result of Nussbaum~\cite{MR1269677},
who showed that a modification of this formula, with an absolute value
enclosing each term 
$\frac{\partial \phi_i}{\partial x_j}(s,x)$ for $i\neq j$,
holds for a non-necessarily order-preserving flow.
Nussbaum's approach uses the fact that the Thompson metric originates from a Finsler structure to determine the local contraction rate. However,
this leads to different assumptions
(see Assumption H3.1 in~\cite{MR1269677}).
In particular, as in the discrete case (see Remark~\ref{rem-er}), the method of~\cite{MR1269677} 
requires some form of geodesic convexity assumption,
which can be dispensed with if the flow is assumed
to be order-preserving.
For instance, only the special case of Corollary~\ref{coroNus} in which the
domain $\cU$ is geodesically convex can be recovered by the method
of
~\cite{MR1269677}. 

\subsection{Standard Riccati operator}\label{subsec-std}
One major application of the above analysis is the Riccati operator, arising from the Linear Quadratic (LQ) control problem. Let $\E$ be a real or complex Hilbert space with inner product $\<\cdot,\cdot>$. The set of bounded linear operators on $\E$ is denoted by $\End(\E)$. For $A\in \End(\E)$, let $A'$ denote the adjoint of $A$. The set of symmetric bounded linear operators is denoted by $\ssym(\E)$. A symmetric bounded linear operator $A$ is positive semi-definite if $\<x,Ax>\geq 0$ for all $x\in \E$. Following~\cite{MR2338433}, let $\pP$(resp. $\pP_0$) be the set of positive semi-definite(resp. positive semi-definite invertible) bounded symmetric linear operators of $\E$. Then $\pP$ is a convex closed pointed cone with interior $\pP_0$(Lemma 9.2 and Proposition 9.5~\cite{MR2188393}) and induces the \textit{Loewner} order '$\leq$' on $\ssym(\E)$:
$$P\leq Q \Longleftrightarrow Q-P\in \pP. $$ Then we may define the Thompson metric on $\pP_0$. This is of course a special case of the definition in Section~\ref{subsec-def-th}. Note that equipped with the operator norm, the cone $\pP$ is normal. Therefore the metric space $(\pP_0;d_T)$ is complete(Lemma 5.1~\cite{MR2338433}). 

Consider the Riccati differential equation defined on $\ssym(\E)$:
\begin{align}\label{edfdg}
\dot P(t)=\phi(t,P):=A(t)'P(t)+P(t)A(t)+D(t)-P(t)\Sigma(t)P(t).
\end{align}
where $A:\R\rightarrow \End(\E)$, $D:\R\rightarrow \ssym(\E)$, $\Sigma:\R\rightarrow \ssym(\E)$ are assumed to be continuous and bounded applications.
The flow associated to~(\ref{edfdg}) is naturally order-preserving on $\pP_0$ by considering the LQ control problem. One may also verify it using Proposition~\ref{croissance}.
The least contraction rate of the flow on $\pP_0$ is the best constant $\alpha$ such that for all $P_1,P_2\in \pP_0$ and $s\geq 0$,
\begin{align}\label{leastcontflowcase}
d_T(M_s^t(P_1),M_s^t(P_2))\leq e^ {-\alpha (t-s)}d_T(P_1,P_2), \enspace \forall s \leq t< t_{\pP_0}(s,P_1)\wedge t_{\pP_0}(s,P_2) .
\end{align}
An immediate consequence of Theorem~\ref{nonexp} is:
\begin{theo}\label{tsdfgg}The least contraction rate defined as in~\eqref{leastcontflowcase} satisfies:
\begin{align}\label{a-ekdo}
\alpha=\sup\{\beta\in \R: P\Sigma(t)P+D(t)\geq \beta P,\enspace\forall t\geq 0,P\in\pP_0\}\end{align}
\end{theo}
\begin{rem}\label{rtesfg}
 Even if in the statement of Theorem~\ref{tsdfgg} we do not require $\Sigma$ and $D$ to be positive semi-definite, the set of the supremum of which is taken in~\eqref{a-ekdo} is easily seen to be empty as soon as $\Sigma$ or $D$ are not positive semi-definite. Hence, the finiteness of the constant $\alpha$ in Theorem~\ref{tsdfgg} does require $\Sigma$ and $D$ to be positive semi-definite and then we must have $\alpha\geq 0$. This shows a dichotomy: either the flow is non-expansive, or it is not uniformly Lipschitz.
\end{rem}
\begin{coro}[Theorem 8.5~\cite{MR2338433}]\label{dgfghfgx}
 We suppose that $D(t),\Sigma(t)\in \pP$ for all $t\geq 0$. Then the least contraction rate is given by:
$$
\alpha=2\inf_{t\geq 0}\sqrt{m((\Sigma(t) ^{1/2}D(t)\Sigma(t) ^{1/2})/I)}
$$
\end{coro}
\begin{proof}
The best contraction rate is given by:
$$\alpha=\sup\{\beta\geq 0: P\Sigma(t)P+D(t)\geq \beta P,\enspace\forall t\geq 0,P\in\pP_0\}$$
Consider all $P=\lambda I$, then
$$
\alpha \leq \sup\{\beta\geq 0: \lambda^2\Sigma(t)\geq \beta \lambda I-D(t),\enspace\forall t\geq 0,\lambda > 0\}
$$
If $\Sigma(t)\in \pP$ is not invertible, then $m(\Sigma/I)=0$. Thus
$$
\alpha \leq \sup\{\beta\geq 0: 0\geq \beta \lambda +m(-D(t)/I),\enspace\forall t\geq 0,\lambda > 0\}=0
$$
Now suppose that $\Sigma(t)\in \pP_0$, $\forall t\geq 0$. In that case, $P\Sigma(t) P+D(t)\geq \beta P$ if and only if
$$
\begin{array}{l}
\Sigma(t)^{\frac{1}{2}}P\Sigma(t)P\Sigma(t)^{\frac{1}{2}}+\Sigma(t)^{\frac{1}{2}}D(t)\Sigma(t)^{\frac{1}{2}}-\beta \Sigma(t)^{\frac{1}{2}}P\Sigma(t)^{\frac{1}{2}}\\
=(\Sigma(t)^{\frac{1}{2}}P\Sigma(t)^{\frac{1}{2}})^2-\beta\Sigma(t)^{\frac{1}{2}}P\Sigma(t)^{\frac{1}{2}}+\Sigma(t)^{\frac{1}{2}}D(t)\Sigma(t)^{\frac{1}{2}}\\
=(\Sigma(t)^{\frac{1}{2}}P\Sigma(t)^{\frac{1}{2}}-\frac{\beta}{2}I)^2+\Sigma(t)^{\frac{1}{2}}D(t)\Sigma(t)^{\frac{1}{2}}-\frac{\beta^2}{4}I
\geq 0
\end{array}
$$
Therefore,
$$
\begin{array}{ll}
\alpha&=\sup\{\beta\geq 0:\beta \leq 2\sqrt {m((\Sigma(t)^ {1/2}D(t)\Sigma(t)^ {1/2})/I)},\enspace\forall t\geq 0\}\\
&=2\displaystyle\inf_{t\geq 0}\sqrt{m((\Sigma(t) ^{1/2}D(t)\Sigma(t) ^{1/2})/I)}.
\end{array}
$$
\end{proof}
The above theorem was proved by Lawson and Lim in~\cite{MR2338433}, Theorem 8.5, using a Birkhoff contraction formula of the fractional transformation on symmetric cones. Their approach requires the coefficients $\Sigma(t)$ and $D(t)$ to be positive semi-definite. By Remark~\ref{rtesfg}, this condition is also necessary to the existence of a global contraction rate. However, a local contraction may occur even the coefficients are not positive semi-definite. 

We now consider the Riccati equation with constant coefficients $(A,D,\Sigma)$. The common fixed point of the flow $M_t$ for all $t$ must satisfy the algebraic Riccati equation(ARE) equation:
$$
A'P+PA+D-P\Sigma P=0
$$
If $\Sigma,D\in \pP_0$, by Corollary~\ref{dgfghfgx} and the completeness of the metric space $(\pP_0;d_T)$ we know that the solution of ARE exists and is unique. We next give sufficient conditions for the existence of solutions of ARE even when $\Sigma$ is not positive semi-definite. 
Below is a direct consequence of Theorem~\ref{theo-nondp2}.
\begin{coro}\label{rsdght}
 Let $P_0\in\pP_0$ and $\alpha \in \R$. The following are equivalent:  
\begin{itemize}
\item[(a)]
For all $P_1,P_2\in (0,P_0)$,
$$
d_T(M_t(P_1,P_2))\leq e^ {-\alpha t}d_T(P_1,P_2), \enspace\forall t< t_{(0,P_0)}(P_1)\wedge t_{(0,P_0)}(P_2).
$$
\item [(b)] For all $P\in (0,P_0)$,
$$
D+P\Sigma P\geq \alpha P.
$$
\end{itemize}
\end{coro}

In particular, this corollary allows to prove the local contraction property of the Riccati equation~(\ref{edfdg}) when $\Sigma$ is not positive definite. Let $c_A,c_D,m_D,c_\Sigma\in \R$ such that:
$$
A+A'\leq -2c_AI,\enspace m_DI\leq D\leq c_D I, \enspace \Sigma \geq -c_\Sigma I.
$$
The situation considered in the next corollary is motivated by the analysis of a method of reduction of the curse of dimensionality introduced by McEneaney~\cite{curseofdim}. This method applies to a control problem in which one can switch between several linear-quadratic models.

\begin{coro}\label{ptzsdg}
 Suppose that $c_A, c_D>0 ,m_D,c_\Sigma >0$ and
$$c_A^ 2\geq c_D c_\Sigma,\enspace c_\Sigma m_D> (c_A-\sqrt{c_A^ 2-c_Dc_\Sigma})^ 2,$$
then for any $\lambda \in [\frac{c_A-\sqrt{(c_A^ 2-c_Dc_\Sigma)}}{c_\Sigma},\sqrt{\frac{m_D}{c_\Sigma}})$, there is $\alpha\geq (m_D-c_ \Sigma \lambda^ 2)/\lambda$ such that for all $P_1,P_2\in (0,\lambda I]$
$$
d_T(M_t(P_1),d_T(P_2))\leq e^ {-\alpha t}d_T(P_1,P_2),\quad\forall t\geq 0.
$$
In particular, there exists a unique solution $\bar P$ to ARE in $(0,\lambda I]$ and for any $P\in (0,\lambda I]$,
$$
d_T(M_t(P),\bar P)\leq e^ {-\alpha t}d_T(P,\bar P),\enspace\forall t\geq 0.
$$ 
\end{coro}

\begin{proof}Let any $\lambda \in [\frac{c_A-\sqrt{(c_A^ 2-c_Dc_\Sigma)}}{c_\Sigma},\sqrt{\frac{m_D}{c_\Sigma}})$.
Since
$$
\begin{array}{l}
\phi(\lambda I)=\lambda (A+A')+D-\lambda^ 2\Sigma \leq (-2\lambda c_A+c_D+\lambda^2c_\Sigma)I\leq 0.
\end{array}
$$
we deduce that the closed set $(0,\lambda I]$ is invariant by the Riccati flow. 
It is not difficult to show that given $\lambda_0\in (\lambda,\sqrt{\frac{m_D}{c_\Sigma}})$ there is $\alpha\geq (m_D-c_\Sigma \lambda_0^2)/ \lambda_0$ such that
$$
D+P\Sigma P\geq \alpha P,\enspace\forall P\in (0,\lambda_0 I).
$$
Indeed, note that a sufficient condition would be:
$$
m_D I-c_\Sigma P^ 2\geq \alpha P, \quad\forall P\in (0,\lambda_0 I)
$$
which is equivalent to:
$$
m_D-c_\Sigma \lambda_0^ 2\geq \alpha \lambda_0.
$$
By Corollary~\ref{rsdght}, for any $P_1,P_2\in (0,\lambda I]\subset (0,\lambda_0 I)$,
$$
d_T(M_t(P_1),M_t(P_2))\leq e^ {-\alpha t}d_T(P_1,P_2), \enspace\forall t\geq 0
$$
Since the metric space $((0,\lambda I];d_T)$ is complete, we deduce that there is a unique fixed point $\bar P\in(0,\lambda I]$ and all solutions with initial value in $(0,\lambda I]$ converge exponentially to $\bar P$ with rate $\alpha$.
\end{proof}

Another interesting case is when $\Sigma \in \pP$ not invertible. In that case, Corollary~\ref{dgfghfgx} tells that the least contraction rate on $\pP_0$ is 0. However, using Corollary~\ref{rsdght} we can say something more about the asymptotic behavior of the trajectories.
\begin{coro}\label{erfgghfg}
 Suppose that $c_A>0, m_D>0$ and $c_\Sigma=0$. Then for any $\lambda\geq \frac{c_D}{c_A}$, there is $\alpha>0$ such that the flow is $\alpha$-contractive on the set $(0,\lambda I]$. In particular, the existence and uniqueness of solution $\bar P\in \pP_0$ to ARE is insured and for any $P\in \pP_0$, 
$$
d_T(M_t(P),\bar P)\leq e^ {-\alpha t}d_T(P,\bar P),\enspace\forall t\geq 0.
$$
where $\alpha=\min(m(I/P),\frac{c_A}{c_D})m_D$.
\end{coro} 
We leave the proof to the reader, which is similar to the one of Corollary~\ref{ptzsdg}.

\if

Next we consider the discrete Riccati operator $F: \pP_0\rightarrow \pP_0$:
\begin{align}\label{rzdfgh}
F(P)=Q+A(B+P^{-1})^ {-1}A'
\end{align}
where $Q,B\in \pP$ and $A\in \End(\E)$, the Lipschitz constant of $F$ is defined by:
\begin{align}\label{ruapfg}
\mylip{F}:=\displaystyle\sup_{P_1,P_2\in \pP_0} \frac{d_T(F(P_1),F(P_2))}{d_T(P_1,P_2)}.
\end{align}
By Proposition~\ref{pazeff}, it is immediate that the operator $F$ is order-preserving on $\pP_0$. Now applying Proposition~\ref{ptthsd}, we get
\begin{align}\label{perrsd}
\mylip{F}=\inf\{\alpha \geq 0: DF(P)\cdot P\leq \alpha F(P),\forall P\in \pP_0\}.
\end{align}
\begin{coro}[Compare with~\cite{MR2338433}]\label{pzsffgg}
Let $\mylip{F}$ be the Lipschitz constant of $F$. 
\begin{itemize}
\item[(1)] If $B\in \pP_0$, then,
 $$\mylip{F}=\frac{\nu}{(1+\sqrt{1+\nu})^ 2}$$
where $\nu=M(AB^ {-1}A'/Q)$.
\item[(2)] If $A$ is invertible and $B$ is not invertible, then
$$\mylip{F}=1.$$
\end{itemize}
\end{coro}
\begin{proof}
(1) By~(\ref{perrsd}), $\mylip{F}$ is the least $\alpha\geq 0$ such that:
$$
A(B+P^ {-1})^ {-1}P^ {-1}(B+P^ {-1})^ {-1}A'\leq \alpha(Q+A(B+P^ {-1})^ {-1}A'),\enspace\forall P\in \pP_0.
$$
Equivalently, $1-\mylip{F}$ is the largest $\eta\leq 1$ such that:
$$
-A(B+P^{-1})^ {-1}B(B+P^ {-1})^ {-1}A'\leq (1-\eta)Q-\eta A(B+P^ {-1})^ {-1}A',\enspace\forall P\in\pP_0.
$$
Let $V=B^ {1/2}(B+P^ {-1})^ {-1}B^ {1/2}$. Then after a simple manipulation we get that $\mylip{F}$ is the largest $\eta<1$ such that:
$$
(1-\eta)Q-AB^ {-\frac{1}{2}}(\eta V-V^ 2)B^ {-\frac{1}{2}}A'\geq 0,\enspace\forall V\in (0,I).
$$
Whence we have $\eta\in [0,1]$. Now since 
$$\eta V-V^ 2\leq \frac{\eta^ 2}{4}I,\enspace\forall V\in (0,I),$$
and when $V=\frac{\eta}{2}I\in (0,I)$,
$$
\eta V-V^2=\frac{\eta^2}{4}I,
$$
$1-\mylip{F}$ is the largest $\eta\in [0,1]$ such that:
\begin{align}\label{afdrr}
(1-\eta)Q\geq \frac{\eta^ 2}{4}AB^ {-1}A',
\end{align}
therefore
$$
\mylip{F}=\frac{\nu}{(1+\sqrt{1+\nu})^ 2}
$$
where $\nu=M(AB^ {-1}A'/Q)$.

(2), If $A$ is invertible, by~(\ref{perrsd}), $\mylip{F}$ is the least $\alpha\geq 0$ such that:
$$
(B+P^ {-1})^ {-1}P^ {-1}(B+P^ {-1})^ {-1}\leq \alpha(A^ {-1}QA'^ {-1}+(B+P^ {-1})^ {-1}),\enspace\forall P\in \pP_0.
$$
Let $\gamma=M(A^ {-1}QA'^ {-1}/I)$, put $P^ {-1}=\lambda I\in \pP_0$, then
\begin{align}
\mylip{F}&\leq \inf\{\alpha\geq 0: \lambda(B+\lambda I)^ {-2}\leq \alpha (\gamma I+(B+\lambda I)^ {-1}),\forall \lambda>0\}\\
&\leq \inf\{\alpha\geq 0: \alpha \gamma^ 2B^ 2+\alpha(2\gamma\lambda+1)B+(\alpha \gamma\lambda^ 2-(1-\alpha)\lambda)I\geq 0,\forall \lambda>0 \}
\end{align}
Therefore if $B$ is not invertible, $\mylip{F}=1$.
If $B$ is invertible, as in the case (1), we have that $1-\mylip{F}$ is the largest $\eta\in [0,1]$ such that~(\ref{afdrr}) is satisfied,
which implies $\mylip{F}=1$ if $Q$ is not invertible.
\end{proof}
\begin{rem}
When $A$ is invertible, the operator~(\ref{rzdfgh}) can be seen as a linear fractional transformation given by symplectic Hamiltonian operators. The contraction coefficient of this operator with respect to Thompson's metric was already calculated by a Birkhoff formula both in finite~\cite{liveraniw} and infinite dimensional case (Theorem 7.2~\cite{MR2338433}). When $A$ is not invertible, the finite dimension case was treated in ~\cite{MR2399829} (estimation (5.16) in Theorem 5.3). Our results show that the last estimation is tight and generalize it to the infinite dimensional case.
\end{rem}

\fi

\section{Application to stochastic Riccati differential equations}\label{sec-sto}

In the sequel, denote by $\symdim{n}$ the set of $n$-dimensional real symmetric matrices. Observe that equipped with the canonical inner product 
$$
\<A,B>=\Tr(AB), \enspace \forall A,B\in \symn,
$$
$\symn$ is a real Hilbert space. The subset of all positive semi-definite matrices $\symnpos\subset \symn$ forms a closed pointed convex cone. All the above results can then be applied here considering $\symnpos$ as $\C$ and the set of all positive definite matrices $\symndefpos$ as $\C_0$. Note that here $\C^*=\C$. We shall use the notation $\geq$ 
(and $\gg$) for the (strict) Loewner order,
 and $d_T$ for the Thompson metric induced by $\symnpos$ (see Section~\ref{subsec-def-th}).

\label{sertcgg}
\subsection{Stochastic LQ problem and GRDE}
Consider the following stochastic linear quadratic optimal control problem:
\begin{displaymath}\label{stoLQ}
\begin{array}{l}
v(s,y)=\displaystyle\min_{u(\cdot)} \E \int_{s}^T[x(t)'Q(t)x(t)+2u'(t)L(t)x(t)+u(t)'R(t)u(t)]dt+\E [x(T)'Gx(T)]\\
\textrm{s.t.}\quad \left \{\begin{array}{l}
dx(t)=(A(t)x(t)+B(t)u(t))dt+(C(t)x(t)+D(t)u(t))dW(t),\quad \forall t\in [s,T],
\\
x(s)=y.
\end{array}\right.
\end{array}
\end{displaymath}
where the functions appearing above satisfy:
$$
\left\{\begin{array}{l}
A(\cdot),C(\cdot)\in L^{\infty}\cap C^0(0,T;\R^{n\times n}),\quad B(\cdot),D(\cdot),L(\cdot)\in L^{\infty}\cap C^{0}(0,T;\R^{n\times k}),\\ Q(\cdot)\in L^{\infty}\cap C^{0}(0,T;\symn),\quad R(\cdot)\in L^{\infty}\cap C^{0}(0,T;\symdim{k}).\end{array}\right.$$ 
Here $W$ is a standard Brownian motion defined on a complete probability space. We refer the reader to~\cite{YongZhoubook99} Chapter 6, for the precise definition of this control problem.
In~\cite{YongZhoubook99}, the above functions are only assumed to be bounded. In our case, the continuity is necessary to apply the previous results.

The above stochastic LQ control problem over the time interval $[s,T]$ is solvable,
i.e., admits an optimal control for all $y\in \R^n$ if the solution of the following constrained differential matrix equation exists:
\begin{align}\label{GRDE2}
\left\{\begin{array}{l}
\dot P+A'P+PA+C'PC+Q=\\ \qquad \quad(PB+C'PD+L')(R+D'PD)^{-1}(B'P+D'PC+L), \quad \enspace t\in [s,T]\\
P(T)=G\\
R(t)+D(t)'P(t)D(t)\gg 0,\quad  \enspace t\in [s,T]
\end{array}\right.
\end{align}
which we refer to as generalized Riccati differential equation (GRDE). 
In that case, the value function of the optimal control problem is given by \begin{align}\label{a-valu}v(s,y)=y'P(s)y.\end{align}

\if
The partial order induced by the cone $\symnpos$ is the Loewner order:
$$
\begin{array}{l}
A\geq B \Leftrightarrow A-B \in \symnpos\\
A\gg B \Leftrightarrow A-B \in \symndefpos
\end{array}
$$
The Thompson part metric is well defined on $\symndefpos$ and can be calculated by:
$$d_T(A,B)=\log \max\{\rho(A^{-1}B),\rho(B^{-1}A)\}$$
where $\rho(\cdot)$ denote the spectral radius.
\fi
\subsection{GRDE with semi-definite weighting matrices}
The solvability of the GRDE~(\ref{GRDE2}) with indefinite matrix coefficients has been treated by Chen, Moore, Rami, and Zhou in~\cite{MR1819815}. In order to apply our previous results, we only consider the case:
\begin{align}\label{positivecondition}\left(\begin{array}{ll}Q(t)& L(t)' \\L(t)&R(t)\end{array}\right)\geq 0,\enspace \ker{R(t)}\cap \ker{D(t)}=\{0\},\enspace \forall t\in [0,T],\end{align}
so that the function
\begin{align}\label{functionphi}
\begin{array}{ll}
\phi(t,P)=&PA+A'P+C'PC+Q-\\&(B'P+D'PC+L)'(R+D'PD)^{-1}(B'P+D'PC+L)
\end{array}
\end{align}
is well defined on $[0,+\infty)\times\symndefpos$ and satisfies the assumptions made at the beginning of section~\ref{se12}.
We are going to apply the preceding results to show the monotonicity and the non-expansiveness of the GRDE differential equation defined on $\symndefpos$:
\begin{align}\label{GRDE}
\left\{\begin{array}{l}
\dot P=\phi(t,P), \quad\\
P(0)=G\\
\end{array}\right.
\end{align}
\begin{prop}\label{stononexp}
Assume that~\eqref{positivecondition} holds. Then the flow associated to~(\ref{GRDE}) is order-preserving and non-expansive on $\symndefpos$. 
\end{prop}

 Proposition~\ref{stononexp} could be derived by exploiting the relation between the solution of the Riccati equation and the value function of the stochastic control problem (see~\eqref{a-valu}). Here we choose to prove it from the infinitesimal characterizations of Proposition~\ref{croissance} and Theorem~\ref{nonexp}.

\begin{proof}
By Proposition~\ref{croissance}, if suffices to prove that for any $P\in \symndefpos$, any $Q,Z\in \symnpos$ such that $\<Q,Z>=0$:
$$\<Q,D\phi_t(P) Z>\geq 0.$$
Indeed, 
$$
\begin{array}{ll}
D\phi_t(P) Z=&ZA(t)+A(t)'Z+C(t)'ZC(t)-(B(t)'Z+D(t)'ZC(t))'N_t(P)\\
&-N_t(P)'(B(t)'Z+D(t)'ZC(t))+N_t(P)'D(t)'ZD(t)N_t(P)
\end{array}
$$
where $N_t(P)=(R(t)+D(t)'PD(t))^{-1}(B(t)'P+D(t)'PC(t)+L(t))$.
Remark that if $Q,Z\in \symnpos$ and $\<Q,Z>=0$ then $QZ=0$. Therefore,
$$
\begin{array}{ll}
\<Q,D\phi_t(P) Z>&=\langle Q,C(t)'ZC(t)-C(t)'ZD(t)N_t(P)-N_t(P)'D(t)'ZC(t)
\\ &\quad +N_t(P)'D(t)'ZD(t)N_t(P)\rangle\\
\qquad\qquad\qquad\enspace&=\langle Q, \big (C(t)-D(t)N_t(P)\big)'Z \big (C(t)-D(t)N_t(P)\big)\rangle \geq 0.
\end{array}
$$
Now for non-expansiveness, by Theorem~\ref{nonexp} it remains to verify that for any $P\in \symndefpos$ and any $t\in [0,T]$,$$
D\phi_t(P) P-\phi(t,P)\leq 0.$$ 
Indeed,
\begin{align}\label{aftehig}
\begin{array}{ll}
&D\phi_t(P) P-\phi(t,P)\\
&=-Q(t)+N_t(P)'L(t)+L'(t)N_t(P)-N_t(P)'R(t)N_t(P)\\
&=H_t(P)' \left(\begin{array}{ll}-Q(t)& -L(t)' \\-L(t)&-R(t)\end{array}\right)H_t(P)
\leq 0
\end{array}
\end{align}
where $H_t(P)'=\left(
\begin{array}{ll}
I&-N_t(P)'
\end{array}
\right)$.
\end{proof}
\begin{rem}
A fundamental discrepancy with the standard Riccati equation is that the flow of the generalized Riccati equation is not a global contraction. This is because that there is no $\alpha>0$ such that the condition
$$
D\phi_t(P) P-\phi(t,P)\leq -\alpha P,\enspace \forall P\in\symndefpos,
$$
which by Theorem~\ref{nonexp} is necessary to the global contraction property of the flow, is satisfied. However, we shall see in the next section that a local contraction property does hold.
\end{rem}

\subsection{Asymptotic behavior of GRDE}
We are going to investigate the behavior of the GRDE flow as time horizon goes to infinity. All the matrices $A,B,C,D,L,Q,R$ are assumed to be constant. First we show a local contraction property under the condition \begin{align}\label{arsdf}\left(\begin{array}{ll}Q& L' \\L&R\end{array}\right)\gg 0.\end{align} More precisely, 

\begin{theo}\label{localcontraction}
Assume that \eqref{arsdf} holds. Let $\cU\subset \symndefpos$ be an open set such that $\lambda \cU\subset \cU$ for all $\lambda \in(0,1]$. Assume that there is $P_0\in \symndefpos$ such that $\cU\subset (0, P_0]$
and let $\alpha=m(Q-L'R^{-1}L/P_0)$,
then for all $P_1,P_2\in \cU$,
$$
d_T(M_t(P_1),M_t(P_2))\leq e^{-\alpha t} d_T(P_1,P_2),\enspace  0\leq t<t_\cU(P_1)\wedge t_\cU(P_2)
$$
\end{theo}

\begin{proof}
 By applying Theorem~\ref{theo-nondp2}, we need to prove
$$
D\phi(P)P-\phi(P)\leq -\alpha P,\quad \forall P\in \cU
$$
Indeed, 
for all $P\in \cU$,
$$
Q-\alpha P-L'R^{-1}L\geq Q-\alpha P_0-L'R^{-1}L\geq 0.
$$
Besides,
the previous calculus yields
\begin{align}\label{fffaftehig}
\begin{array}{ll}
&D\phi(P) P-\phi(P)+\alpha P\\
&=H(P)' \left(\begin{array}{ll}-Q+\alpha P & -L' \\-L&-R\end{array}\right)H(P)
\end{array}
\end{align}
where $H(P)'=\left(
\begin{array}{ll}
I&-N(P)'
\end{array}
\right)$ and $N(P)=(R+D'PD)^{-1}(B'P+D'PC+L)$. 
By Schur's complement lemma, we get
$$
D\phi(P) P-\phi(P)\leq -\alpha P,\quad \forall P\in \cU.
$$
\end{proof}

The fixed point of the GRDE flow associated to~(\ref{GRDE}), if it exists, satisfies the so-called general algebraic Riccati equation (GARE):
\begin{align}\label{GARE}
\left\{\begin{array}{l}
\phi(P)=0.\\
R+D'PD\gg0
\end{array}\right.
\end{align}
where $\phi(P):=A'P+PA+C'PC+Q-(B'P+D'PC+L)'(R+D'PD)^{-1}(B'P+D'PC+L)$.
The existence of solutions of GARE and the asymptotic behavior of the GRDE flow have been studied in~\cite{MR1819815} and \cite{MR1778371}. 
The authors assumed the following mean-square stabilizability condition:
\begin{defi}[Definition 4.1~\cite{MR1778371}]
 The system of matrices $(A,B,C,D)$ is said to be mean-square stabilizable if there exists a control law of feedback form
$$
u(t)=Kx(t),
$$
where $K$ is a constant matrix, such that for every initial $(t_0,x_0)$, the closed loop system
$$
\left\{\begin{array}{l}
dx(t)=(A+BK)x(t)dt+(C+DK)x(t)dW(t)\\
x(0)=x_0
\end{array}\right.
$$
satisfies$$
\lim_{t\rightarrow+\infty}\E[x(t)'x(t)]=0
$$
\end{defi}
Under the mean-square stabilizability assumption, they established a necessary and sufficient condition for the existence of a solution. To make a comparison, let us first quote their theorem:
\begin{theo}[Theorem 4.1~\cite{MR1819815}]
 Under the mean-square stabilizability assumption, there exists a solution of the GARE ~(\ref{GARE}) if and only if there exists $P_0\in\symn$ such that
$$
\phi(P_0)\geq0, R+D'P_0D\gg0
$$
Moreover, for any such $P_0$, the solution $P(t)$ of~(\ref{GRDE}) with initial condition $P(0)=P_0$ converges to a solution to the GARE as $t\rightarrow \infty$. 
\end{theo}
It follows directly from the above theorem that under the mean-square stabilizability assumption, if \eqref{arsdf} is true, then there must be a solution to the GARE~(\ref{GARE}). We next show a necessary and sufficient condition for the existence of a stable solution without the mean-square stabilizability assumption.


\begin{theo}\label{gareNEW} Assume that the condition~\eqref{arsdf} holds. Then, the GARE admits a solution $\bar P\in\symndefpos$ if and only if there exists $P_0\in \symndefpos$ such that:
\begin{align}\label{dfezq}
\phi(P_0)\leq 0.
\end{align}
In that case, for any $P\in \symndefpos$:
$$
d_T(M_t(P),\bar P)\leq e^{-\alpha t} d_T(P,\bar P),\enspace \forall t\geq 0,
$$
where $$\alpha\geq \frac{1-e^{-d_T(P,\bar P)}}{d_T(P,\bar P)}m((Q-L'R^{-1}L)/\bar P)>0.$$ In particular, the solution is unique in $\symndefpos$.
\end{theo}
\begin{proof}
 If $\bar P\in\symndefpos$ is a solution of the GARE, then~(\ref{dfezq}) is satisfied by considering $P_0=\bar P$. Conversely, note that if $\phi(P_0)\leq 0$ for some $P_0\in \symndefpos$, then $(0,P_0]$ is an invariant set. Consider the open set $\cU=(0,P_0+I)$. 
By Theorem~\ref{localcontraction}, there is $\alpha>0$ such that for all $P_1$, $P_2\in (0,P_0]\subset \cU$, we have:
$$
d_T(M_t(P_1),M_t(P_2))\leq e^{-\alpha t} d_T(P_1,P_2), \enspace\forall 0\leq t\leq t_{\cU}(P_1)\wedge t_{\cU}(P_2),
$$
Since $[0,P_0]\subset \cU$ is invariant, we have that $t_{\cU}(P_1),t_{\cU}(P_2)=+\infty$.
Thus the flow $M_t$ is contractive in the complete metric space $((0,P_0],d_T)$. There must be a unique fixed point $\bar P\in (0,P_0]$ such that $\phi(\bar P)=0$. 
Next, assuming the existence of a solution $\bar P\in \C_0$ to the GARE, we apply Theorem~\ref{tzsdfg} to obtain the rate of convergence. A basic calculus yields:
$$
\begin{array}{l}
\lambda^{-1}\phi(\lambda \bar P)\\=(B'\bar P+D'\bar P\tilde C)'((R+D'\bar PD)^{-1}-(\lambda^{-1}R+D'\bar PD)^{-1})(B'\bar P+D'\bar P\tilde C)+(\lambda^{-1}-1)\tilde Q
\end{array}
$$
where $\tilde C=C-DR^{-1}L$ and $\tilde Q=Q-L'R^{-1}L$. 
Therefore, if $\lambda\geq 1$, then
\begin{align}\label{tdfsd}
\lambda^{-1}\phi(\lambda \bar P)\leq (\lambda^{-1}-1)\tilde Q
\end{align}
and \begin{align}\label{tfds}
\lambda \phi(\lambda^{-1}\bar P)\geq (\lambda-1)\tilde Q.
\end{align}
Now for any $P\in \symndefpos \neq \bar P$, let $\mu=e^{d_T(P,\bar P)}$ and $\alpha=\frac{1-\mu^{-1}}{\ln\mu}m(\tilde Q/\bar P)>0$. Then~
$$
(\lambda^{-1}-1)\tilde Q\leq -\alpha (\ln \lambda) \bar P, (\lambda-1) \tilde Q \geq \alpha (\ln \lambda) \bar P, \enspace\forall \lambda \in(1,\mu)
$$
and \eqref{tdfsd} and~\eqref{tfds} lead to:
$$
\begin{array}{ll}
\alpha \ln(\lambda) \bar P\leq \lambda \phi(\lambda^ {-1}\bar P), \enspace \alpha \ln(\lambda) \bar P\leq -\lambda^ {-1}\phi(\lambda P),\enspace\forall \lambda \in(1,\mu).
\end{array}
$$
Thus,
\begin{align}\label{esdfert}
0<\alpha \leq \inf_{\mu^{-1}< \lambda < \mu} m\big((-(\lambda \ln \lambda)^{-1} \phi(\lambda \bar{P}) )/ \bar{P}\big) \enspace .
\end{align}
By virtue of~\eqref{esdfert} and Theorem~\ref{tzsdfg}, we have
$$
d_T(M_t(P),\bar P)\leq e^ {-\alpha t}d_T(P,\bar P),\enspace \forall t\geq 0.
$$
\end{proof}


\subsection{Discrete Generalized Riccati operator}\label{subsec-gen-disc}
The linear quadratic stochastic control problem has a discrete time analogue~\cite{981058}, which leads to the generalized discrete Riccati operator $F:\symn \rightarrow \symn$:
\begin{align}\label{atrrtfg}
F(P)=A'PA+C'PC+Q-(B'PA+D'PC)'(R+B'PB+D'PD)^{-1}(B'PA+D'PC)
\end{align}
where $A,C\in\R^{n\times n}$, $B,D\in \R^{n\times m}$ and $Q, R\in \symn$. We assume that $Q\gg0$ and $R\gg 0$. Then by applying the Schur complement condition for positive definiteness, one can prove that $F$ sends $\symndefpos$ to itself.
Note that when $C=D=0$, 
we recover the standard Riccati operator:
\begin{align}\label{alig-stan-ric}
 T(P)=A'PA+Q-A'PB(R+B'PB)^{-1}B'PA.
\end{align}
The object of this section is to get the Lipschitz constant of $F$ on $\symndefpos$ (see~\eqref{ruapfg-nocited}). First we show that this operator is order-preserving on $\symndefpos$.
\begin{prop}
 The operator $F$ is order-preserving on~$\symndefpos$.
\end{prop}
\begin{proof}
Let any $P\in \symndefpos$ and $Z\in \symnpos$. A simple calculus show that:
 $$
DF(P)\cdot Z=(A-BN)'Z(A-BN)+(C-DN)'Z(C-DN)\geq 0
$$
where $N=(R+B'PB+D'PD)^{-1}(B'PA+D'PC)$. By Proposition~\ref{croissance}, $F$ is order-preserving on $\symndefpos$.
\end{proof}
Next we apply Proposition~\ref{ptthsd} to get:
\begin{align}\label{perrsd}
\mylip{F;\symndefpos}=\inf\{\alpha \geq 0: DF(P)\cdot P\leq \alpha F(P),\forall P\in \symndefpos\}.
\end{align}
The following two lemmas will be useful.
\begin{lemma}\label{ltzrfgcv}
 Let $\bigl( \begin{smallmatrix}
  B\\ D
\end{smallmatrix} \bigr)
=\bigl( \begin{smallmatrix}
  \bar B\\ \bar D
\end{smallmatrix} \bigr)\begin{smallmatrix}
  W
\end{smallmatrix}$ be a rank factorization (so that the last two factors have maximal column and row rank, respectively). Then the operator $F$ defined in~(\ref{atrrtfg}) satisfies:
\begin{align}\label{artdfggg}
F(P)=A'PA+C'PC+Q-(\bar B'PA+\bar D'PC)'(\bar R+\bar B'P\bar B+\bar D'P\bar D)^{-1}(\bar B'PA+\bar D'PC)
\end{align}
where $\bar R=(WR^{-1}W')^{-1}$.
\end{lemma}
\begin{proof}
To simplify the notation, denote $X(P)=\bar B P\bar B+\bar D'P\bar D$. Notice that since the matrix $\bigl( \begin{smallmatrix}
 \bar B\\\bar D
\end{smallmatrix} \bigr)$ is of full column rank, $X(P)$ is invertible for all $P\in \symndefpos$.
It follows from~(\ref{atrrtfg}) that:
$$
F(P)=A'PA+C'PC+Q-(\bar B'PA+\bar D'PC)'W(R+W'X(P)W)^{-1}W'(\bar B'PA+\bar D'PC)
$$
Now appealing to the Woodbury matrix identity, we obtain:
\begin{align}
 \begin{array}{ll}
  W(R+X(P)W)^{-1}W'&=W(R^{-1}-R^{-1}W'(X(P)^{-1}+WR^{-1}W')^{-1}WR^{-1})W'\\
&=WR^{-1}W'-WR^{-1}W'(X(P)^{-1}+WR^{-1}W')^{-1}WR^{-1}W'\\
&=((WR^{-1}W')^{-1}+X(P))^{-1}
 \end{array}
\end{align}
from which we get~(\ref{artdfggg}).
\end{proof}
\begin{lemma}\label{lgty}
Let $\delta\geq 2$, then
\begin{align}\label{aitufg}
 X-X(R+X)^{-1}(\delta R+X)(R+X)^{-1}X\leq \frac{R}{4(\delta -1)},\enspace \forall X\in \symnpos \end{align}
\end{lemma}
\begin{proof}
 Let any $X\in \symnpos$. Since $X$ commutes with $I$, we have that:
$$
\begin{array}{l}
X-X(I+X)^{-1}(\delta I+X)(I+X)^{-1}X \\
=(I+X)^{-1}(X(I+X)^2-X^2(\delta I+X))(I+X)^{-1}\\
=(I+X)^{-1}((2-\delta)X^2+X-\frac{1}{4(\delta-1)}(I+X)^2)(I+X)^{-1}+\frac{1}{4(\delta-1)}I\\
=-(I+X)^{-1}((2\delta-3)X-I)^2(I+X)^{-1}+\frac{1}{4(\delta-1)}I\\
\leq \frac{1}{4(\delta-1)}I.
\end{array}
$$
To obtain~(\ref{aitufg}), it suffices to notice that:
$$
\begin{array}{l}
R^{-\frac{1}{2}}(X-X(R+X)^{-1}(\delta R+X)(R+X)^{-1}X)R^{-\frac{1}{2}}\\
=Y-Y(I+Y)^{-1}(\delta I+Y)(I+Y)^{-1}Y
\end{array}
$$
where $Y=R^{-\frac{1}{2}}XR^{-\frac{1}{2}}$.
\end{proof}

\begin{prop}\label{prop-discrete}
The operator $F$ is non-expansive: $\mylip{F;\symndefpos}\leq 1$. Let 
\[ \left(\begin{array}{l}B\\D\end{array}\right)=
\left(\begin{array}{l}\bar B\\\bar D\end{array}\right)W
\]
be a rank factorization. Then a necessary and sufficient condition to have $\mylip{F;\symndefpos}<1$ is that there is a matrix $S$ such that:
\begin{align}\label{rtsfgfg}
\left(\begin{array}{l}
       A\\C
      \end{array}
\right)=\left(\begin{array}{l}
       \bar B\\\bar D
      \end{array}
\right)S.
\end{align} In that case, $$\mylip{F;\symndefpos}\leq \frac{M(S'\bar RS/Q)}{(1+\sqrt{1+M(S'\bar RS/Q)})^2}<1$$ where $\bar R=(WR^{-1}W')^{-1}$.
\end{prop}
\begin{proof}
Lemma~\ref{ltzrfgcv} implies that it is sufficient to prove the proposition for the case $W=I$, i.e. when $\bigl( \begin{smallmatrix}
  B\\ D
\end{smallmatrix} \bigr)$ is of full column rank.
A simple calculus shows that:
$$
\begin{array}{ll}
DF(P)\cdot P-\alpha F(P)=&(1-\alpha)(A'PA+C'PC)-\alpha Q\\&-(1-\alpha)N(P)'(R+X(P))^{-1}N(P)\\&-N(P)'(R+X(P))^{-1}R(R+X(P))^{-1}N(P)
\end{array}
$$
where $$N(P)=B'PA+D'PC,\quad ~X(P)=B'PB+D'PD.$$
Then it is evident that $\mylip{F;\symndefpos}\leq 1$. 
Now let $S\in \R^{n \times m}$ such that~(\ref{rtsfgfg}) holds. Then $N(P)=X(P)S$, $A'PA+C'PC=S'X(P)S$ and
$$
\begin{array}{ll}
DF(P)\cdot P-\alpha F(P)=&(1-\alpha)S'X(P)S-\alpha Q\\&-(1-\alpha)S'X(P)(R+X(P))^{-1}X(P)S\\&-S'X(P)'(R+X(P))^{-1}R(R+X(P))^{-1}X(P)S
\end{array}
$$
To simplify the notation, let $X:=X(P)$ and $\delta:=\frac{2-\alpha}{1-\alpha}$, then
$$
DF(P)\cdot P-\alpha F(P)=(1-\alpha)S'(X-X(R+X)^{-1}(\delta R+X)(R+X)^{-1}X)S-\alpha Q.
$$
By Lemma~\ref{lgty}:
$$
\begin{array}{l}
X-X(R+X)^{-1}(\delta R+X)(R+X)^{-1}X
\leq \frac{1}{4(\delta-1)}R=\frac{1-\alpha}{4}R, \enspace\forall X\in \symnpos.
\end{array}
$$
Therefore 
$$
(1-\alpha)S'(X-X(R+X)^{-1}(\delta R+X)(R+X)^{-1}X)S\leq \frac{(1-\alpha)^2}{4}S'RS.
$$
Consequently if $\alpha$ is such that: $\frac{4\alpha}{(1-\alpha)^2}=M(S'RS/Q)$, then
$$
DF(P)\cdot P-\alpha F(P)\leq 0, \enspace \forall P\in \symnpos.
$$
Together with~\eqref{perrsd} this shows that $$\mylip{F;\symndefpos}\leq \frac{M(S'RS/Q)}{(1+\sqrt{1+M(S'RS/Q)})^2}.$$
Next we prove the necessity of condition~(\ref{rtsfgfg}). Remember that since the matrix $\bigl( \begin{smallmatrix}
  B\\ D
\end{smallmatrix} \bigr)$ has full rank, $X(P)$ is always invertible for $P\in \symndefpos$.
Besides,
there is $\alpha<1$ such that 
$$DF(P)\cdot P-\alpha F(P)\leq 0,\enspace\forall P\in \symndefpos$$
if and only if for any $P\in \symndefpos$,
$$
\begin{array}{l}
(A'PA+C'PC)-\frac{\alpha}{1-\alpha} Q- N(P)'(R+ X(P))^{-1}(\frac{2-\alpha}{1-\alpha}R+ X(P))( R+ X(P))^{-1} N(P)\leq 0.
\end{array}
$$
That is, for any $P\in \symndefpos$ and $\lambda>0$,
$$
\begin{array}{l}
(A'PA+C'PC)-\frac{\alpha\lambda^{-1}}{1-\alpha} Q- N(P)'(\frac{1}{\lambda}R+ X(P))^{-1}(\frac{2-\alpha}{\lambda(1-\alpha)}R+ X(P))(\frac{1}{\lambda}R+ X(P))^{-1} N(P)\leq 0.
\end{array}
$$
Letting $\lambda$ go to infinity, by continuity, we obtain that:
$$
\begin{array}{l}
(A'PA+C'PC)- N(P)' X(P))^{-1} N(P)\leq 0.
\end{array}
$$
The above expression is the Schur complement of the positive semi-definite matrix
$$
\left(
\begin{array}{ll}
 B'PB+D'PD & B'PA+D'PC\\
A'PB+C'PD & A'PA+C'PC
\end{array}\right)=\left(\begin{array}{l}B'\\A'\end{array}\right)P\left(\begin{array}{ll}B&A\end{array}\right)+\left(\begin{array}{l}D'\\C'\end{array}\right)P\left(\begin{array}{ll}D&C\end{array}\right).
$$
 Therefore for any $x\in\R^n$ there is $u\in \R^m$ such that
\begin{align}
&\< \left(\begin{array}{l}u\\x\end{array}\right),\left(\begin{array}{ll}
 B'PB+D'PD & B'PA+D'PC\\
A'PB+C'PD & A'PA+C'PC
\end{array}\right)\left(\begin{array}{l}u\\x\end{array}\right) >= 0.\end{align}
That is, for any $x\in \R^n$ there is $u\in \R^m$ such that: 
$$
\left(\begin{array}{ll}
B&A\\
D&C
\end{array}\right)\left(
\begin{array}{l}
 u\\x
\end{array}
\right)=0.
$$
This is equivalent to say that there is $S\in \R^{m\times n}$ such that:
$$
\left(\begin{array}{l}
       A\\C
      \end{array}
\right)=\left(\begin{array}{l}
       BS\\DS
      \end{array}
\right).
$$
\end{proof}

The contraction rate of the standard discrete Riccati operator $T: \symndefpos\rightarrow \symndefpos$ can now be recovered as a corollary:
\begin{coro}[Compare with~\cite{MR2399829}]
 The standard Riccati operator $T$ defined in~\eqref{alig-stan-ric} is non-expansive: $\mylip{T;\symndefpos}\leq 1$. A necessary and sufficient condition to have the strict contraction property is that the matrix $B$ is of full row rank. In that case, let $B=\bar B W$ be a rank factorization, then
$$\mylip{T;\symndefpos}\leq \frac{M(S'\bar RS/Q)}{(1+\sqrt{1+M(S'\bar RS/Q)})^2}<1$$ where $S=\bar B^{-1} A$ and $\bar R=(WR^{-1}W')^{-1}$.
\end{coro}

\begin{rem}
Condition~(\ref{rtsfgfg}) leads to a formal argument explaining why strict global contraction cannot be hoped for the GRDE flow. Indeed, we can approximate the continuous-time LQ control problem in Section~\ref{sertcgg} over a small time horizon $\epsilon$ by the following one-step discrete-time stochastic linear quadratic control problem:
$$
\begin{array}{l}
\displaystyle\min_{u\in \R^m} \E(\<x_0,\epsilon Qx_0>+\<u,\epsilon Ru>+\<x_\epsilon,Gx_\epsilon>)\\
~~s.t.~~ x_\epsilon=(I+\epsilon A)x_0+\epsilon Bu+(\sqrt{\epsilon} Cx_0+\sqrt{\epsilon} Du)w
\end{array}
$$
where $w\sim \N(0,1)$. Without loss of generality, we suppose that $\bigl( \begin{smallmatrix}
  B\\ D
\end{smallmatrix} \bigr)$ is of full column rank.
If a strict contraction result was valid for the continuous
time system, we would expect the same to be true for its discrete approximation
if $\epsilon$ is sufficiently small. However, 
the strict global contraction condition requires the existence of $S$ such that:
$$
\left(\begin{array}{l}
       I+\epsilon A\\\ \sqrt \epsilon C
      \end{array}
\right)=\left(\begin{array}{l}
       \epsilon BS\\\sqrt \epsilon DS
      \end{array}
\right),$$
which can not hold for a set of $\epsilon$ converging to 0 if $C$ and $D$ are not zero. 
\end{rem}




\section{Loss of non-expansiveness of the GRDE flow in other invariant Finsler metrics}\label{section-loss}
The standard Riccati flow is known to be a contraction in the standard Riemannian metric~\cite{MR1227540}, and more generally in any invariant Finsler metric (with the same bound on the contraction rate)~\cite{MR2399829}. We next construct an explicit
counter example showing that the Thompson metric is essentially
the only invariant Finsler metric in which the GRDE Riccati flow is
non-expansive.

\subsection{Preliminary results}
We first recall the definition of symmetric gauge functions and of the
associated invariant Finsler metrics on the interior of the cone of positive definite
matrices. Then, we will show some conditions that are necessary for
an order-preserving flow to be non-expansive in a given metric of this kind.


\begin{defi}[Symmetric gauge function]
A symmetric gauge function $\nu:\R^n \rightarrow \R$ is a convex, positively homogeneous of degree 1 function such that for any permutation $\sigma$,
$$
\nu(\lambda_1,\cdots,\lambda_n)=\nu(|\lambda_{\sigma(1)}|,\cdots, |\lambda_{\sigma(n)}|),\enspace\forall\enspace \lambda=(\lambda_1,\cdots,\lambda_n)\in \R^n.
$$ 
\end{defi}
The next lemma collects several useful properties of subdifferentials of symmetric gauge function (see~\cite{MR0274683} for more background on subdifferentials). The straightforward proof is left to the reader.
\begin{lemma}\label{laerdf}
 Let $\nu:\R^n\rightarrow \R$ be a symmetric gauge function. The following properties hold:
\begin{itemize}
 \item [1] For all $\lambda \in \R^n$ and $\mu \in \partial \nu(\lambda)$,
$$\mu_i \lambda_i \geq 0,\quad \forall i=1,\cdots,n.$$
\item [2] For all $\lambda \in \R^n$ and $\mu \in \partial \nu(\lambda)$,
$$
\<\mu,\lambda>=\nu(\lambda).
$$
\item [3] For all $\lambda,\lambda'\in \R^n$ and $\mu \in \partial \nu(\lambda), \mu' \in \partial \nu(\lambda')$,
$$
\<\mu-\mu',\lambda>\geq 0.
$$
\end{itemize}

\end{lemma}

For every symmetric gauge function $\nu$, we define a spectral function $\hat \nu:\symn\rightarrow \R$:
$$
\hat \nu(P)=\nu (\lambda (P)).
$$
where $\lambda (P)$ is the vector of eigenvalues of $P$.
\begin{theo}[\cite{LewisMR1377729}]\label{lewisspectral}
 If $\nu$ is a symmetric gauge function, then $\hat \nu$ is a convex function on $\symn$. Moreover, $Z\in \partial \hat \nu (P)$ if and only if there exists $y\in \partial \nu (\lambda(P))$ such that:
\[
 Z=V\diag(y)V^T,
\]
where $V$ is the unitary matrix such that $P=V\diag(\lambda(P))V^T$.
\end{theo}
Following~\cite{MR2013466}, \cite{MR2399829} and \cite{MR1636922},
we define a metric on $\symndefpos$ as follows,
$$
d_{\nu}(P,Q)=\hat \nu(\log (P^{-1/2}QP^{-1/2})).
$$
It coincides with the Finsler metric obtained by thinking of $\symndefpos$ as a manifold and taking
$$
\|dQ\|_P
=\hat\nu(P^{-\frac{1}{2}}(dQ) P^{-\frac{1}{2}})$$
as the length of an infinitesimal displacement in the tangent space at point $P$. 
This metric
is invariant by the canonical action on the linear group on $\symndefpos$.

We shall consider specially, as in~\cite{MR1636922}, the $p$-norm function:
$$
\nu (\lambda)=\|\lambda\|_p=(\sum_{i=1}^n|\lambda_i|^p)^{1/p},
$$
so that the metric $d_{\nu}$ is the Thompson metric for $p=+\infty$ and the Riemannian metric for $p=2$. 
\begin{lemma}
Let $\nu$ be a symmetric gauge function and $d_{\nu}$ be the associated metric on $\symndefpos$. Let $M:\symnpos\rightarrow \symnpos$ be a differentiable function such that: 
\begin{align} d_\nu(M(P),M(Q))\leq d_\nu(P,Q),\enspace \forall\enspace P,Q\in \symndefpos,\label{nonexnu}
\end{align}then
$$ \hat \nu(M(P)^{-1/2}(D M(P)\cdot Z) M(P)^{-1/2})\leq  \hat \nu (P^{-1/2}ZP^{-1/2}),\enspace \forall  \enspace P\in \symndefpos,Z\in \symn.$$
\end{lemma}
\begin{proof}
Let any $P\in \symndefpos$ and $Z\in \symn$. There exists $\delta>0$ such that for any $0\leq \epsilon\leq \delta$, $P+\epsilon Z\in \symndefpos$.
By (\ref{nonexnu}) and the definition of $d_{\nu}$:
\[
\hat\nu \log (M(P)^{-1/2}M(P+\epsilon Z)M(P)^{-1/2}) \leq \hat\nu \log (P^{-1/2}(P+\epsilon Z)P^{-1/2}).
\]
Divide the two sides by $\epsilon$ and take the limit:
\[
\lim_{\epsilon \rightarrow 0} \frac{\hat\nu \log (M(P)^{-1/2}M(P+\epsilon Z)M(P)^{-1/2})}{\epsilon} \leq \lim_{\epsilon \rightarrow 0} \frac{\hat\nu \log (P^{-1/2}(P+\epsilon Z)P^{-1/2})}{\epsilon}
\]
In view of homogeneity and continuity of the function $\hat \nu$,
\[
\hat \nu (\lim_{\epsilon\rightarrow 0} \frac{\log(M(P)^{-1/2}M(P+\epsilon Z)M(P)^{-1/2})}{\epsilon}) \leq \hat \nu (\lim_{\epsilon\rightarrow 0} \frac{\log P^{-1/2}(P+\epsilon Z)P^{-1/2}}{\epsilon})
\]
The matrix function $\log$ is differentiable at $I$:
\begin{align}
\lim_{\|U\| \rightarrow 0} \frac{\log(I+U)-U}{\|U\|}=0. \label{Dlog}
\end{align}
Hence by chain rule:
\[
\hat \nu (M(P)^{-1/2}(D M(P) \cdot Z)M(P)^{-1/2}) \leq  \hat \nu (P^{-1/2}Z P^{-1/2})
\]
\end{proof}

We consider the following time independent differential equation:
\begin{align}\label{lstsectintial}
\left\{
\begin{array}{l}
\dot x(t)=\Phi(x(t)),
\\ x(s)=x_0.
\end{array}\right.
\end{align} 
where $\Phi$ is differentiable on $\symndefpos$. We assume that the associated flow $M_{\cdot}(\cdot): (0,+\infty)\times \symndefpos\rightarrow \symn $ leaves $\symndefpos$ invariant and is globally defined.

\begin{lemma}
Let $\nu$ be a symmetric gauge function. If there exists $\epsilon>0$ such that for any $0\leq t\leq \epsilon$, 
\begin{align}\label{nonexnc} \hat \nu(M_t(I)^{-1/2}(\cD M_t(I)\cdot Z) M_t(I)^{-1/2})\leq  \hat \nu (Z),\enspace \forall  \enspace Z\in \symn.\end{align}
then
$$
\<\diag(\mu),\cD\Phi(I)\cdot \diag(\lambda)-\diag(\lambda)\Phi(I)>\leq 0,\enspace \forall \lambda \in \R^n,\enspace \mu \in \partial \nu(\lambda).
$$
\end{lemma} 
\begin{proof}
Let any $Z\in \symn$.
For readability, 
denote
$$
P_t:=M_t(I),\enspace
H_t:=P_t^{1/2},\enspace
Q_t:=P_t^{-1/2},\enspace
G_t:=P_t^{-1},$$ and
$$
U_t:=\cD M_t(I)\cdot Z,\enspace J_t:=U_t G_t,\enspace K_t=Q_tJ_tH_t.$$ 
The derivative of $J_t$ with respect to $t$ is:
\begin{align}
\begin{array}{ll}
\dot J_t&=
\dot U_t G_t-U_t G_t \dot P_t G_t       \\
&=(\cD\Phi(P_t)\cdot U_t)G_t-U_t G_t \Phi(P_t)G_t.
      \end{array}
 \label{Jprime}
\end{align}
The derivative of $K_t$ with respect to $t$ is:
\[
\begin{array}{ll}
\dot K_t&=Q_t\dot J_t H_t+\dot Q_t J_t H_t+Q_t J_t \dot H_t\\
&=Q_t \dot J_t H_t -Q_t \dot H_t Q_t J_t H_t+Q_t J_t\dot H_t\\
&=Q_t \dot J_t H_t-Q_t\dot H_t K_t+K_t Q_t \dot H_t
\end{array}
\]
Hence,
$$\dot J_t|_{t=0}=\cD\Phi(I)\cdot Z-Z\Phi(I),$$and
\begin{align}\label{Kt0}
\dot K_t|_{t=0}=\cD\Phi(I)\cdot Z-Z\Phi(I)-(\dot H_t|_{t=0})Z+Z (\dot H_t|_{t=0}).
\end{align}
By Theorem~\ref{lewisspectral}, the right derivative of the function $\hat \nu(K_t)$ with respect to $t$ exists:
\[\begin{array}{ll}
\hat \nu(K_t)'_+=\displaystyle\sup_{y\in \partial \hat \nu(K(t))}\<y,\dot K_t>&=\displaystyle\sup_{\substack{\mu \in \partial \nu(\lambda),V V'=I\\V' K_t V=\diag(\lambda)}}\langle V\diag(\mu)V',\dot K_t\rangle.\\
\end{array}
\]
Since
$$
\hat\nu(K_t)\leq \hat \nu(K_0), \enspace  t\in [0,\delta),
$$
the right derivative at $t=0$ must be negative:
$$
\hat \nu(K_t)'_+|_{t=0}\leq 0.
$$
Namely,
$$
\displaystyle\sup_{\substack{\mu \in \partial \nu(\lambda),V V'=I\\V' Z V=\diag(\lambda)}}\langle V\diag(\mu)V',\cD \Phi(I)\cdot Z-Z\Phi(I)-(\dot H_t|_{t=0})Z+Z (\dot H_t|_{t=0})\rangle\leq 0.
$$
Note that for any unitary matrix $V$ such that $V'ZV=\diag(\lambda)$, we have
$$
\<V\diag(\mu)V',(\dot H_t|_{t=0})Z>=\<V\diag(\mu)V',Z (\dot H_t|_{t=0})>
$$
Hence by taking $Z=\diag(\lambda)$ and $V=I$, we obtain a necessary condition of (\ref{nonexnc}):
$$
\<\diag(\mu),\cD\Phi(I)\cdot \diag(\lambda)-\diag(\lambda)\Phi(I)>\leq 0
$$
for all $\lambda \in \R^n$ and $\mu \in \partial \nu(\lambda)$.
\end{proof}
The above two lemmas lead to the following conclusion:
\begin{prop}
If the flow $M_{\cdot}(\cdot):(0,+\infty)\times \symndefpos\rightarrow \symndefpos $ is non-expansive in the metric $d_{\nu}$, then,
\begin{align}\label{nessarycondition}
\<\diag(\mu),\cD\Phi(I)\cdot \diag(\lambda)-\diag(\lambda)\Phi(I)>\leq 0 
\end{align}
for all $\lambda \in \R^n$ and $\mu \in \partial \nu(\lambda)$.
\end{prop}
\subsection{The counter example}
We finally arrive at the announced counter example:
we give a system of parameters $(A,B,C,D,L,Q,R)$ such that the corresponding $\Phi$ of GRDE does not satisfy the necessary condition~(\ref{nessarycondition}) of non-expansiveness in any Finsler metric other than the Thompson metric.

Recall that 
$$
\Phi(P)=A'P+PA+C'PC+Q-(B'P+D'PC+L)'(R+D'PD)^{-1}(B'P+D'PC+L).
$$
Let $I_n$ denote the $n$-dimensional identity matrix and $e=(e_1,\cdots,e_{n-1})'\in \R^{n-1}$ be a vector.    
The parameters are chosen as follows:
$$
A=I_{n},\enspace
B=\Big(\begin{array}{ll}
(\epsilon-\sqrt{1-\epsilon})I_{n-1} &0\\
-(\sqrt{1-\epsilon}) e'&\epsilon
\end{array}\Big),\enspace
C=\Big(\begin{array}{ll}
(1+\sqrt{1-\epsilon})I_{n-1} & e\\
0&\sqrt{1-\epsilon}
\end{array}\Big),$$and
$$
D=(\sqrt{1-\epsilon})I_{n}
,\enspace
L=0,\enspace
R=\epsilon I_{n},\enspace
Q=\epsilon I_{n}
$$
to make
$$
R+D'D=I_{n},\enspace B'+D'C=I_{n},\enspace
C-D=(\begin{array}{ll}
I_{n-1} & e\\
0&0
\end{array})
$$
An elementary calculus yields
$$
\begin{array}{l}
D\Phi(I)\cdot Z-Z\Phi(I)\\
=A'Z+ZA+C'ZC-(B'Z+D'ZC)'(R+D'D)^{-1}(B'+D'C)\\ \quad-(B'+D'C)(R+D'D)^{-1}(B'Z+D'ZC)\\ \quad +(B'+D'C)(R+D'D)^{-1}D'ZD(R+D'D)^{-1}(B'+D'C)\\
\quad -Z(A'+A+C'C+Q-(B'+D'C)'(R+D'D)^{-1}(B'+D'C))\\
=2Z+C'ZC-(B'Z+D'ZC)'-(B'Z+D'ZC)+D'ZD-Z-ZC'C-Z'Q\\
=Z+(C-D)'Z(C-D)-B'Z-ZB'-ZC'C-ZQ
\end{array}
$$
Now let any $\lambda=(\lambda_1,\cdots,\lambda_{n})\in \R^{n}$ and $\mu \in \partial \nu (\lambda)$. Then
$$
\begin{array}{l}
\<\diag(\mu),\cD\Phi(I)\cdot \diag(\lambda)-\diag(\lambda)\Phi(I)>\\
=-2\epsilon \<\mu,\lambda>+\mu_n (-\lambda_n |e|^2+\sum_{i=1}^{n-1}\lambda_i e_i^2)
\end{array}
$$
Recall that $$\<\mu,\lambda>=\nu(\lambda),\quad \forall \mu\in \partial \nu(\lambda).$$
So if there is any $\lambda\in \R^{n}$ and $\mu \in \partial \nu (\lambda)$ such that
$$
\mu_n(-\lambda_n |e|^2+\sum_{i=1}^{n-1} \lambda_i e_i^2)>0,
$$
then there always exists $\epsilon\in (0,1)$ such that
$$
\<\diag(\mu),\cD\Phi(I)\cdot \diag(\lambda)-\diag(\lambda)\Phi(I)> >0.
$$
Finally we need a lemma to conclude:
\begin{lemma}
 If for all $\lambda \in \R^{n}$, $\mu \in \partial \nu (\lambda)$ and $e\in \R^{n-1}$ we have
$$
\mu_{n}(-\lambda_{n}\|e\|^2+\sum_{i=1}^{n-1} 
\lambda_i e_i^2)\leq 0,
$$
then $$\nu(\lambda_1,\cdots,\lambda_n)=c\max_i |\lambda_i|$$
for some constant $c>0$.
\end{lemma}
\begin{proof}
 First consider $e=e_i$ the $i$-th
standard basis vector of $\R^{n-1}$ for all $i=1,\cdots,n-1$. We see that
$$
\mu_n(-\lambda_n+\lambda_i)\leq 0,\quad \forall  i=1,\cdots,n-1
$$
for all $\lambda=(\lambda_1,\cdots,\lambda_n)\in \R^n$ and $\mu=(\mu_1,\cdots,\mu_n) \in \partial \nu (\lambda)$.
By the symmetric property of $\nu$, this implies actually
\begin{align}\label{afgress}
\mu_j(-\lambda_j+\lambda_i)\leq 0,\quad \forall i,j=1,\cdots,n
\end{align}
Therefore, for any $\lambda\neq 0$ if $\lambda_j=0$ then $\mu_j=0$ for all $\mu\in \partial \nu(\lambda)$. 
Next, let any $i\in \{1,\cdots,n\}$, consider the following set
$$\Lambda_i:=\{\lambda \neq 0: \lambda_1=\lambda_2=\dots=\lambda_i>\lambda_{i+1}\geq \dots \lambda_{n}\geq 0\}.$$
Let any $\lambda\in \Lambda_i$ and $\mu \in \partial \nu(\lambda)$. By Property 1 in Lemma~\ref{laerdf}, $\mu \geq 0$. Using~\eqref{afgress}, we know that:
$$\mu_j\leq 0,\quad\forall j=i+1,\dots,n.$$
Hence,
$$
\mu_j=0,\quad \forall j=i+1,\dots,n.
$$
Now let any $\lambda^1,\lambda^2\in \Lambda_i$ and $\mu^1\in \partial \nu(\lambda^1)$, $\mu^2\in \partial \nu(\lambda^2)$. By Property 3 in Lemma~\ref{laerdf},
$$
\<\mu^1-\mu^2,\lambda^1>\geq 0.
$$
It follows that
$$
\sum_{j=1}^i \mu^1_j\geq \sum_{j=1}^i \mu^2_j.
$$
We deduce that $\sum_{j=1}^i \mu^1_j=\sum_{j=1}^i \mu^2_j$. Hence there is a constant $c_i\geq 0$ such that
$$
\nu(\lambda)=\<\mu,\lambda>=\sum_{i=1}^j \mu_j \lambda_j=\lambda_1\sum_{i=1}^j \mu_j=c_i \lambda_1, \quad \forall \lambda \in \Lambda_i,\mu \in \partial \nu(\lambda).
$$ 
It remains to prove that $c_i=c_1$ for all $i=1,\dots,n$. To see this, again we use Property 3 in Lemma~\ref{laerdf}. First consider $\lambda=(1,\dots,1,0,\dots,0)\in \Lambda_i$ and any $\mu \in \partial \nu(\lambda)$, then
$$
\<\mu-(c_1,0,\dots,0)',\lambda>=\sum_{j=1}^i \mu_j -c_1=c_i-c_1 \geq 0.
$$
On the other hand, for all $\lambda^1\in \Lambda_1 $
$$
\<(c_1,0,\dots,0)'-\mu,\lambda^1>=(c_1-\mu_1)\lambda^1_1-\sum_{j=2}^i \mu_j \lambda_j^1\geq 0
$$
This implies $\displaystyle c_1=\sum_{j=1}^i \mu_j=c_i$ for all $i=1,\dots,n$.
\end{proof}

The proof of Theorem~\ref{theo-final} is now complete.

%


\section*{Acknowledgments}
The authors thank Shanjian Tang for having raised the problem of the contraction rate of the generalized Riccati flow and also for several suggestions all along the course of the present work.

\bibliographystyle{alpha}
\bibliography{biblio}

\begin{thebibliography}{RCMZ01b}

\bibitem[ACS00]{MR1636922}
E.~Andruchow, G.~Corach, and D.~Stojanoff.
\newblock Geometrical significance of {L}\"owner-{H}einz inequality.
\newblock {\em Proc. Amer. Math. Soc.}, 128(4):1031--1037, 2000.

\bibitem[Bha03]{MR2013466}
Rajendra Bhatia.
\newblock On the exponential metric increasing property.
\newblock {\em Linear Algebra Appl.}, 375:211--220, 2003.

\bibitem[Bir57]{birkhoff57}
Garrett Birkhoff.
\newblock Extensions of {J}entzsch's theorem.
\newblock {\em Trans. Amer. Math. Soc.}, 85:219--227, 1957.

\bibitem[Bon69]{bony}
J.-M. Bony.
\newblock Principe du maximum, in\'egalit\'e de harnack et unicit\'e du
  probl\`eme de cauchy pour les op\'erateurs elliptiques d\'eg\'en\'er\'es.
\newblock {\em Annales de l'institut Fourier}, 19(1):277--304, 1969.

\bibitem[Bou93]{MR1227540}
Philippe Bougerol.
\newblock Kalman filtering with random coefficients and contractions.
\newblock {\em SIAM J. Control Optim.}, 31(4):942--959, 1993.

\bibitem[Bre70]{MR0257511}
Ha{\"{\i}}m Brezis.
\newblock On a characterization of flow-invariant sets.
\newblock {\em Comm. Pure Appl. Math.}, 23:261--263, 1970.

\bibitem[Cla75]{MR0367131}
Frank~H. Clarke.
\newblock Generalized gradients and applications.
\newblock {\em Trans. Amer. Math. Soc.}, 205:247--262, 1975.

\bibitem[CLSW98]{Clarke:1998:NAC:274798}
F.~H. Clarke, Yu.~S. Ledyaev, R.~J. Stern, and P.~R. Wolenski.
\newblock {\em Nonsmooth analysis and control theory}.
\newblock Springer-Verlag New York, Inc., Secaucus, NJ, USA, 1998.

\bibitem[CLZ98]{ChengLiZhou98}
Shuping Chen, Xunjing Li, and Xun~Yu Zhou.
\newblock Stochastic linear quadratic regulators with indefinite control weight
  costs.
\newblock {\em SIAM J. Control Optim.}, 36(5):1685--1702 (electronic), 1998.

\bibitem[Lew96]{LewisMR1377729}
A.~S. Lewis.
\newblock Convex analysis on the {H}ermitian matrices.
\newblock {\em SIAM J. Optim.}, 6(1):164--177, 1996.

\bibitem[LL06]{MR2188393}
Jimmie Lawson and Yongdo Lim.
\newblock The symplectic semigroup and {R}iccati differential equations.
\newblock {\em J. Dyn. Control Syst.}, 12(1):49--77, 2006.

\bibitem[LL07]{MR2338433}
Jimmie Lawson and Yongdo Lim.
\newblock A {B}irkhoff contraction formula with applications to {R}iccati
  equations.
\newblock {\em SIAM J. Control Optim.}, 46(3):930--951 (electronic), 2007.

\bibitem[LL08]{MR2399829}
Hosoo Lee and Yongdo Lim.
\newblock Invariant metrics, contractions and nonlinear matrix equations.
\newblock {\em Nonlinearity}, 21(4):857--878, 2008.

\bibitem[LW94]{liveraniw}
Carlangelo Liverani and Maciej~P. Wojtkowski.
\newblock Generalization of the {H}ilbert metric to the space of positive
  definite matrices.
\newblock {\em Pacific J. Math.}, 166(2):339--355, 1994.

\bibitem[Mar73]{MR0318991}
R.~H. Martin, Jr.
\newblock Differential equations on closed subsets of a {B}anach space.
\newblock {\em Trans. Amer. Math. Soc.}, 179:399--414, 1973.

\bibitem[McE07]{curseofdim}
W.~M. McEneaney.
\newblock A curse-of-dimensionality-free numerical method for solution of
  certain {HJB} {PDE}s.
\newblock {\em SIAM J. Control Optim.}, 46(4):1239--1276, 2007.

\bibitem[Nus88]{nussbaum88}
R.~D. Nussbaum.
\newblock Hilbert's projective metric and iterated nonlinear maps.
\newblock {\em Mem. Amer. Math. Soc.}, 75(391):iv+137, 1988.

\bibitem[Nus94]{MR1269677}
Roger~D. Nussbaum.
\newblock Finsler structures for the part metric and {H}ilbert's projective
  metric and applications to ordinary differential equations.
\newblock {\em Differential Integral Equations}, 7(5-6):1649--1707, 1994.

\bibitem[RCMZ01a]{aitrami01}
Mustapha~Ait Rami, Xi~Chen, John~B. Moore, and Xun~Yu Zhou.
\newblock Solvability and asymptotic behavior of generalized {R}iccati
  equations arising in indefinite stochastic {LQ} controls.
\newblock {\em IEEE Trans. Automat. Control}, 46(3):428--440, 2001.

\bibitem[RCMZ01b]{MR1819815}
Mustapha~Ait Rami, Xi~Chen, John~B. Moore, and Xun~Yu Zhou.
\newblock Solvability and asymptotic behavior of generalized {R}iccati
  equations arising in indefinite stochastic {LQ} controls.
\newblock {\em IEEE Trans. Automat. Control}, 46(3):428--440, 2001.

\bibitem[RCZ01]{981058}
M.A. Rami, X.~Chen, and X.Y. Zhou.
\newblock Discrete-time indefinite lq control with state and control dependent
  noises.
\newblock In {\em Decision and Control, 2001. Proceedings of the 40th IEEE
  Conference on}, volume~2, pages 1249 --1250 vol.2, 2001.

\bibitem[Red72]{MR0303024}
R.~M. Redheffer.
\newblock The theorems of {B}ony and {B}rezis on flow-invariant sets.
\newblock {\em Amer. Math. Monthly}, 79:740--747, 1972.

\bibitem[Roc70]{MR0274683}
R.~Tyrrell Rockafellar.
\newblock {\em Convex analysis}.
\newblock Princeton Mathematical Series, No. 28. Princeton University Press,
  Princeton, N.J., 1970.

\bibitem[RW75]{RedhefferandWalterMR0470401}
R.~M. Redheffer and W.~Walter.
\newblock Flow-invariant sets and differential inequalities in normed spaces.
\newblock {\em Applicable Anal.}, 5(2):149--161, 1975.

\bibitem[RZ00]{MR1778371}
Mustapha~Ait Rami and Xun~Yu Zhou.
\newblock Linear matrix inequalities, {R}iccati equations, and indefinite
  stochastic linear quadratic controls.
\newblock {\em IEEE Trans. Automat. Control}, 45(6):1131--1143, 2000.

\bibitem[Tho63]{ThompsonMR0149237}
A.~C. Thompson.
\newblock On certain contraction mappings in a partially ordered vector space.
\newblock {\em Proc. Amer. Math. soc.}, 14:438--443, 1963.

\bibitem[YZ99]{YongZhoubook99}
Jiongmin Yong and Xun~Yu Zhou.
\newblock {\em Stochastic controls}, volume~43 of {\em Applications of
  Mathematics (New York)}.
\newblock Springer-Verlag, New York, 1999.
\newblock Hamiltonian systems and HJB equations.

\end{thebibliography}

\end{document}